\documentclass[journal]{IEEEtran}
\usepackage[dvips]{graphicx}
\usepackage[english]{babel}
\usepackage[latin1]{inputenc}
\usepackage{amssymb}
\usepackage{url}
\usepackage{cite}
\usepackage{stfloats}
\usepackage{color,soul}
\newcommand{\resethlcolor}{\sethlcolor{yellow}}
\resethlcolor
\definecolor{lightblue}{rgb}{.90,.95,1}
\definecolor{bluegreen}{rgb}{0.2,0.5,0.2}

\begin{document}
\newtheorem{theorem}{Theorem}
\newtheorem{lemma}{Lemma}
\newtheorem{corollary}{Corollary}
\newtheorem{proposition}{Proposition}
\newtheorem{definition}{Definition}
\newcommand{\vv}[1]{\mathbf{#1}}
\newcommand{\E}{{\bf E}}
\newcommand{\Prob}{{\bf P}}

\title{ Asymptotically Optimal Change Point Detection for Composite Hypothesis in State Space Models}

\author{Cheng-Der~Fuh 
\thanks{Cheng-Der Fuh is with the Fanhai International School of Finance, Fudan University, China 
(email: cdfu@fudan.edu.cn). }}

%\thanks{The research of Cheng-Der Fuh is partially supported by MOST 105-2410-H-008-025-MY2, and MOST 106-2118-M-008-002-MY2. }}
\date{November 20, 2018}

\maketitle

\begin{abstract}

This paper investigates change point detection in state space models, in which the pre-change distribution $f^{\theta_0}$ is given, while the poster distribution $f^{\theta}$ after change is unknown. The problem is to raise an alarm as soon as possible after the distribution changes from $f^{\theta_0}$ to $f^{\theta}$, under a restriction on the false alarms. 
We investigate theoretical properties of a weighted 
Shiryayev-Roberts-Pollak (SRP) change point detection rule in state
space models. By making use of a Markov chain representation for the likelihood function, exponential embedding of the 
induced Markovian transition operator, nonlinear Markov renewal theory, and sequential hypothesis testing theory for Markov random walks, we show that the weighted SRP procedure is second-order asymptotically optimal. To this end, we derive an asymptotic approximation for the expected stopping time of such a
stopping scheme when the change time $\omega = 1$. To illustrate our  method we apply the results to two types of state space models: general state Markov chains and linear state space models.
\end{abstract}
\centerline{{\bf Index Terms}}~

Asymptotic optimality, change point detection, first passage time,
iterated random functions system, nonlinear Markov renewal theory, sequential analysis, Shiryayev-Roberts-Pollak procedure.

\section{Introduction}\label{Intro}
\def\theequation{1.\arabic{equation}}
\setcounter{equation}{0}

A prototypical problem of detecting abrupt changes can be found in
instruction detection in distributed computer networks. Large scale attacks, denial of service attacks, occur at unknown points in time and need to be detected at the early stages
by observing abrupt changes in the computer network traffic.
Further applications are in, for example, biomedical signal processing, industrial quality control, segmentation of signals, financial engineering, edge detection in images, and the diagnosis of faults in the elements of computer communication networks.
The reader is referred to Lai \cite{Lai1995, Lai1998} and Tartakovsky et al. \cite{TNBbook14} for a comprehensive summary in this area. 
A standard formulation of the change point detection problem is that there is a sequence of observations whose distribution changes at some unknown time $\omega$, and the goal is to detect this change as soon as possible under false alarm constraints.

When the observations $Y_n$ are independent with a common density
function $f^{\theta_0}$ for $n < \omega$ and with another common density function $f^{\theta}$ for $n \geq \omega$, where $\omega$ is unknown and both $\theta_0$ and $\theta$ are given, there are two standard formulations for the optimum tradeoff problem. The first is a minimax formulation proposed by Lorden \cite{Lorden1971},
in which he shows that subject to the
``average run length'' (ARL) constraint, Page's CUSUM procedure asymptotically minimizes the ``worst case'' detection delay. The second is a Bayesian formulation, proposed by
Shiryayev \cite{Shiryayev1963, Shiryayev1978}, in which the change point has a geometric prior distribution on,
and the goal is to minimize the expected delay subject to an upper bound on false alarm probability.
He uses optimal stopping theory to show that the Bayes rule triggers an alarm as soon as the posterior probability that a change has occurred exceeds some fixed level. Roberts \cite{Roberts1966} considers the non-Bayesian setting, and studies by simulation the average run length of this rule, and finds 
it to be very good. Pollak \cite{Pollak1985} shows that the (modified) 
Shiryayev-Roberts rule is asymptotically minimax.
When $\theta$ is unknown, Pollak and Siegmund \cite{Pollak1975} extends Shiryayev's work in a non-Bayesian setting,
and calculates the expected value of a weighted likelihood ratio test as well as the average run lengths of a CUSUM rule. Then Pollak \cite{Pollak1987} provides average run lengths of the weighted Shiryayev-Roberts change point detection rule.

Regarding change point detection rules in dynamic systems beyond independent assumption. In the case of using CUSUM type change point detection rules, Bansal and Papantoni-Kazakos \cite{Bansal1986} extends Lorden's asymptotic theory to the case where $Y_j$ are stationary ergodic sequences, under the condition that 
$\{Y_j, j < \omega\}$ (before the change point) and $\{Y_j, j \geq
\omega\}$ (after the change point) are independent, 
and proves the asymptotic optimality of the CUSUM algorithm. Further extensions to general stochastic sequences $Y_n$ were obtained by Lai \cite{Lai1995, Lai1998}, 
and Tartakovsky and Veeravalli \cite{Tartakovsky2005}. 
When both $\theta_0$ and $\theta$ are given,
Fuh \cite{Fuh2003} proves that the CUSUM scheme is asymptotically optimal, in the sense of Lorden \cite{Lorden1971},
in hidden Markov models. In the domain of Shiryayev-Roberts type change point detection rules, Yakir \cite{Yakir1994} generalizes the
result to a finite state Markov chain, while Bojdecki \cite{Bojdecki1979} studies a different loss function and applies optimal stopping theory to find the Bayes rule. 
Tartakovsky \cite{Tartakovsky-IEEEIT2017}  considers a sequential Bayesian changepoint detection problem for a general stochastic model.
Fuh \cite{Fuh2004b} investigates the Shiryayev-Roberts-Pollak 
(SRP) change point detection rule in hidden Markov models, in which he proves the asymptotic minimax property and derives an asymptotic approximation for the average run lengths when $\omega =1$. 
Fuh and Tartakovsky \cite{FuhTa2018} considers asymptotic Bayesian change point detection in hidden Markov models.

It is noted that many practical problems for change point detection
are beyond independent assumption. Some useful
class of such models are $AR$ models, $ARMA$ models, and linear state space models, cf. Tartakovsky et al. \cite{TNBbook14}. Along this line, in this paper, we study change point detection in state 
space models. A prototypical {\it state space model} 
can be formulated as follows: for $n=1,2,\ldots,$ define
\begin{eqnarray}\label{ssm1}
 Y_n = G_{\theta}(X_n,\varepsilon_n), ~~~{\rm and}~~~ X_n = F_\theta(X_{n-1},\eta_n),
\end{eqnarray}
where $Y_n$ is the observed value, $X_n$ is a $d$-dimensional vector representing an unobservable
state, and $(\varepsilon_n, \eta_n)$ are independent random vectors representing random disturbances
and having a common density function $\phi_{\theta}$. 
Furthermore, we assume $\{\varepsilon_n, n \geq 0\}$ and $\{\eta_n, n \geq 0\}$ are independent. Here the system dynamics are given by the second equation in (\ref{ssm1}). Note that the state vectors $X_n$ are not directly observable and the observations are $Y_n$ which
are related to $X_n$ and measurement error $\varepsilon_n$ in the first equation of (\ref{ssm1}). 

Specifically, a simple linear state space model given by MacGregor and Harris (1990) 
to study the problem of monitoring process means with the sample means $Y_n$:
\begin{eqnarray}\label{mh}
Y_n = X_n + \varepsilon_n,~~~{\rm and}~~~X_n-\mu=\alpha(X_{n-1} - \mu) + \eta_n,
\end{eqnarray}
where $\varepsilon_n$ and $\eta_n$ are independent normal random variables with zero means, and $var(\varepsilon_n)
=\sigma_\varepsilon^2$ and $var(\eta_n)=\sigma_\eta^2$. Here $\theta=(\mu,\alpha,\sigma_\varepsilon^2,\sigma_\eta^2)$ with $|\alpha|<1$, and the
target value of the production process is $\mu=\mu^*$. If we are interested primarily in shifts in the overall
mean and treat $\alpha,\sigma_\varepsilon^2$ and $\sigma_\eta^2$ as unknown nuisance parameters, then we have an incomplete base-line
information and can apply the change point detection rule, described in Section \ref{secAO}, to this case.

In this paper, we will study the change point problem that the pre-change is given while the after change is unknown. It is reasonable to assume that the pre-change distribution is known, because in most practical applications, a large amount of data generated by the pre-change distribution is available to the observer who may use this data to obtain an accurate approximation of the pre-change distribution. However, estimating or even modelling the post-change distribution is often impractical as we may not know a priori what kind of change will happen. We seek to design a change point detection algorithm that allows us to quickly detect the change, under false alarm constraints, and with suitable knowledge of the post-change distribution. To this end, the primary goal of this paper 
is to investigate theoretical properties of a weighted 
Shiryayev-Roberts-Pollak (SRP) change point detection rule in state space models. 

There are two main contributions in this study.
First, we consider a state space model (\ref{ssm1}) in which the underlying state space
is neither finite nor compact, and includes (finite state) hidden Markov models,
linear state space model, and $AR/ARMA$ models as special cases.
Second, the parameter of the distribution after change is assumed to be unknown for practical applications.

The remainder of the paper is organized as follows.
Our main results are in Section \ref{secAO}, in which we derive a second-order asymptotic approximation
for the expected stopping scheme when $\omega = 1$, not worst case, and prove the weighted SRP rule is second-order asymptotically optimal under a false-alarm constraint. In Section \ref{secEA} we illustrate our method by considering two interesting examples: general state Markov chains and linear state space models.
Section \ref{secLR} presents the pre-required methods used in the proofs of our results.
We first give a Markov chain representation of the likelihood ratio, and then study exponential embedding for the induced Markovian transition kernel in state space
models. Based on a nonlinear Markov renewal theory, we
characterize the constant term of the second order approximation in Section \ref{sec5}.
The proofs are given in Sections \ref{sec6}, \ref{sec7} and Appendix, respectively.  

\def\theequation{2.\arabic{equation}}
\setcounter{equation}{0}
\section{Asymptotic optimality of the weighted SRP detection procedure}\label{secAO}

In this section, we define a state space model as a parameterized Markov random walk, in which the underlying environmental Markov
chain can be viewed as a latent variable. To be more precise, for each $\theta \in \Theta \subset {\bf R}$, the unknown
parameter, let ${\bf X}= \{X_n, n \geq 0 \}$ be a Markov
chain on a general state space ${\cal X}$, with transition
probability kernel $P^{\theta}(x,\cdot)= P^{\theta}\{X_1 \in \cdot|X_0=x\}$ and
stationary probability $\pi_{\theta}(\cdot)$. Suppose that a random
sequence $\{Y_n\}_{n=0}^{\infty},$ taking values in ${\bf R}^d$,
is adjoined to the chain such that $\{(X_n,Y_n), n \geq 0\}$ is a
Markov chain on ${\cal X} \times {\bf R}^d$
satisfying $P^{\theta} \{ X_1 \in A | X_0=x,Y_0=y \} = P^{\theta} \{ X_1 \in A | X_0=x \}$
for $A \in {\cal B}({\cal X})$, the $\sigma$-algebra of ${\cal X}$. And conditioning on the full
${\bf X}$ sequence, we have
\begin{eqnarray} \label{ssm}
&~& P^{\theta} \{Y_{n+1} \in  B | X_0,X_1,\ldots;Y_0,Y_1,\ldots, Y_n \} \\
&=& P^{\theta} \{Y_{n+1} \in  B | X_{n+1} \}= P^{\theta} (X_{n+1}: B )~~~a.s. \nonumber
\end{eqnarray}
for each $n$ and $B \in {\cal B}({\bf R}^d),$ the Borel $\sigma$-algebra of ${\bf R}^d$. 
Furthermore, we assume the existence of
a transition probability density $p_{\theta}(x,y)$ for the Markov chain $\{X_n, n \geq 0\}$
with respect to a $\sigma$-finite measure $m$ on ${\cal X}$ such that
\begin{eqnarray}\label{ssmden}
&~& P^{\theta} \{X_1 \in A, Y_{1} \in  B | X_0=x  \} \\
&=& \int_{x' \in A} \int_{y \in B} p_{\theta}(x,x') f(y;\theta|x')  Q(d y) m(dx'), \nonumber
\end{eqnarray}
for $B \in {\cal B}({\bf R}^d)$.
Here $f(Y_k;\theta|X_k)$ is the conditional probability
density of $Y_k$ given $X_k$, with respect to a
$\sigma$-finite measure $Q$ on ${\bf R}^d$.
We also assume that the Markov chain $\{(X_n,Y_n), n \geq 0\}$ has
a stationary probability with probability density function
$\pi_\theta(x)f(\cdot;\theta|x)$ with respect to $m \times Q$.
For convenience of notation, we will use $\pi(x)$ for $\pi_{\theta}(x)$, $p(x,x')$ for
$p_{\theta}(x,x')$, and $f(Y_k|X_k)$ for $f(Y_k;\theta|X_k)$, respectively,
here and in the sequel. We give a formal definition as follows.
\begin{definition}
$\{Y_n,n \geq 0\}$ is called a state space model if there
is an unobserved Markov chain $\{X_n,n \geq 0\}$ such that the process
$\{(X_n,Y_n),n \geq 0\}$ satisfies {\rm (\ref{ssm})}.
\end{definition}

Note that the general state space model defined in Definition 1 includes (\ref{ssm1}), $ARMA$ models, $(G)ARCH$ models and
stochastic volatility models. cf. Fan and Yao \cite{Fan2003} and Fuh \cite{Fuh2006}.

To formulate the change point detection problem,
let $Y_1,\ldots,Y_{\omega -1}$ be a sequence of random variables from the state space model $\{Y_n,n \geq 1\}$ with distribution $P^{\theta_0}$, and let
$Y_{\omega},Y_{\omega+1},\ldots$ be a sequence of random variables from the state space 
model $\{Y_n,n \geq 1\}$ with distribution $P^{\theta}$ at some unknown time $\omega$.
The parameter of pre-change $\theta_0 \in \Theta \subset {\bf R}$ is given;
while the parameter of after change $\theta \in J=(a,b) \subset \Theta$ is unknown.
Moreover, we assume $\theta_0 < a < b < \infty$.
We shall use ${P}_{\omega}$ to denote such a probability measure
(with change time $\omega$) and use ${P}_{\infty}$ to denote the case
$\omega = \infty$ (no change point). Denote $E_{\omega}$ as the corresponding expectation under $P_{\omega}$.
The objectives are to raise an alarm as soon as possible after the change and to avoid false
alarms. A sequential detection scheme $N$ is a stopping time on the sequence of observations $\{Y_n,n \geq 1\}$.
A false alarm is raised whenever the detection is declared before the change occurs. A good detection procedure should minimize the number of post change observations, 
provided that there is no false alarm, while
the rate of false alarms should be low. Hence, the stopping time $N$ should satisfy $\{N \geq \omega \}$ but,
at the same time, keep $N-\omega $ small. Specifically we will find a stopping time $N$ to minimize
\begin{eqnarray}\label{aoe}
\sup_{1 \leq k < \infty} \sup_{\theta \in J} E_k^{\theta}(N-k|N \geq k )
\end{eqnarray}
subject to
\begin{eqnarray}\label{ao1}
E_{\infty}^{\theta_0} N  \geq \gamma,
\end{eqnarray}
for some specified (large) constant $\gamma$. A detection scheme is called second-order asymptotically optimal,
if it minimizes (\ref{aoe}), within an $O(1)$ order, among all stopping rules that satisfy (\ref{ao1}), where $O(1)$ converges to a constant as $\gamma \to \infty$. 
%Note that here we use the global probability of false alarm constraint (\ref{ao1}). 
%Asymptotic optimality of change
%point dectection for stationary sequences, based on
%local probability of false alarm constraint, can be found in Lai (1995, 1998).  
%Asymptotic optimality based on the average run length constraint, $E_\infty N \geq T$ for some specified (large) constant $T$, is still an open problem.

When both $\theta_0 \in \Theta$ and $\theta \in J \subset \Theta$ are given,
the Shiryayev-Roberts-Pollak change point detection scheme
in state space models can be described as follows.
Let $Y_1,\ldots,Y_{n}$ be a sequence of random variables
from the state space model $\{Y_n,n \geq 1\}$, denote
\begin{eqnarray}\label{lr}
&& LR_n (\theta) := \frac{p_n(Y_{1},\ldots,Y_n;\theta)}
{p_n(Y_{1},\ldots,Y_n;\theta_0)}  \\
&&:= \frac{\int_{x_0 \in {\cal X}, \ldots, x_n \in {\cal X}} \pi_{\theta}(x_0)
   \prod_{l=1}^n  p_{\theta}(x_{l-1},x_l) f(Y_l;\theta|x_l)} 
{\int_{x_0 \in {\cal X}, \ldots, x_n \in {\cal X}}  \pi_{\theta_0}(x_0)
  \prod_{l=1}^n  p_{\theta_0}(x_{l-1},x_l) f(Y_l;\theta_0|x_l)}   \nonumber \\
  &&~~~\times \frac{m(dx_n)\cdots m(dx_0)}{m(dx_n)\cdots m(dx_0)} \nonumber
\end{eqnarray}
as the likelihood ratio.
For $0 \leq k \leq n$, denote the detection scheme as
\begin{eqnarray}\label{lrka}
&& LR_n^k (\theta) := \frac{p_n(Y_k,Y_{k+1},\ldots,Y_n;\theta)}
 {p_n(Y_k,Y_{k+1},\ldots,Y_n;\theta_0)}    \\
&&:= \frac{\int_{x_k \in {\cal X}, \ldots, x_n \in {\cal X}}
 \prod_{l=k}^n  p_{\theta}(x_{l-1},x_l) f(Y_l;\theta|x_l)} 
{\int_{x_k \in {\cal X}, \ldots, x_n \in {\cal X}}  \prod_{l=k}^n  p_{\theta_0}(x_{l-1},x_l)
f(Y_l;\theta_0|x_l)}. \nonumber \\
 &&~~~\times \frac{m(dx_n)\cdots m(dx_k)}{m(dx_n)\cdots m(dx_k)} \nonumber
\end{eqnarray}
Given an approximate threshold $B > 0$ and setting $b=\log B$, define the Shiryayev-Roberts scheme as
\begin{eqnarray} \label{sr}
 N_b(\theta) &:=& \inf\{n: \sum_{k=0}^n LR_n^k (\theta) \geq B \} \\
 &=& \inf\{n: \log \sum_{k=0}^n LR_n^k (\theta) \geq b \}. \nonumber
\end{eqnarray}

A simple modification of (\ref{sr}) was given in Pollak \cite{Pollak1985} by adding a randomization 
on the initial $LR_n^0(\theta)$.
This is the celebrated Shiryayev-Roberts-Pollak (SRP) change point detection scheme. 
An extension to finite state hidden Markov models can be found in Fuh \cite{Fuh2004b}.

When $\theta_0 \in \Theta$ is given and $\theta \in J$ is unknown, we apply a similar idea as that in Pollak and Siegmund \cite{Pollak1975},
and Pollak \cite{Pollak1987} for independent observations, extending (\ref{sr}) to have a weight function of $LR_n^k (\theta)$
\begin{eqnarray}\label{mlrk}
LR_n^k(F) := \int_{\theta \in J} LR_n^k (\theta) dF(\theta),
\end{eqnarray}
where $F$ is a probability measure on $J$ with $F(\{\theta_0\})=0$. Given an approximate
threshold $B > 0$ and setting $b=\log B$, define
\begin{eqnarray} \label{msr}
 N_b(F) &:=& \inf\{n: \sum_{k=0}^n LR_n^k (F) \geq B \} \\
 &=& \inf\{n: \log \sum_{k=0}^n LR_n^k (F) \geq b \}. \nonumber
\end{eqnarray}
Then (\ref{msr}) is the weighted SRP change point detection rule in state space models. 
A formal definition will be given in Section 5, in which we will show that the SRP scheme is
an ``equalizer rule'' in the sense that
$E_k ( N_{b}(\theta) - k+1| N_{b}(\theta) \geq k-1) = E_1 N_{b}(\theta),$ for all $k>1$.
%\cdf{double check the stationary property.}

REMARK 1. Note that the weighted SRP change point detection rule (\ref{msr}) involves
two mixture components. One is an integration over the unknown parameter $\theta$
with respect to a prior distribution. The other is an integration over unknown states in the 
state space models, which is related to the non-linear filtering problem. In practice,
it is usually difficult to carry out the computation of $LR_n^k(F)$ in (\ref{mlrk}). A natural substitution is to replace it by
$LR_n^k(\hat{\theta}_{l,k})$ with $\hat{\theta}_{l,k}$ is an estimator of $\theta$ based on $Y_k,\ldots,Y_{l-1}$,
then apply Markov chain Monte Carlo method, in particular particle filtering algorithm, 
to approximate the change point detection rule (\ref{msr}). Theoretical justification and 
empirical study of this change point detection rule are interesting tasks for further investigation.

To derive asymptotic approximation of the average run length, and to prove asymptotic optimality
of the weighted SRP rule in state space models, the following
condition C will be assumed throughout this paper. Before that, we need some definitions first.

A Markov chain $\{X_n,n \geq 0\}$ on a state space ${\cal X}$ is called $V$-uniformly ergodic if there exists a measurable function $V: {\cal X} \rightarrow [1,\infty)$, with
$\int V(x) m(dx) < \infty$, and
\begin{eqnarray}\label{3.11}
& ~&\lim_{n \rightarrow \infty} \sup_{x \in {\cal X}} \bigg\{\frac{\big|E[h(X_n)|X_0=x]
 - \int h(x')m(dx')\big|}{V(x)}: \nonumber \\
& ~&~~~~~~~~~~~~~~~~|h| \leq V \bigg\}=0.
\end{eqnarray}

A Markov chain $\{X_n,n \geq 0\}$ is called Harris recurrent if there exist a recurrent set ${\cal R} \in {\cal B}({\cal X})$, a probability measure $\varphi$ on ${\cal R}$
and an integer $n_0$ such that $P\{ X_n \in {\cal R}~{\rm for~some}~n \geq n_0 |X_0=x\} =1$,
for all $x \in {\cal X}$, and there exists $\lambda > 0$ such that
\begin{eqnarray}\label{min}
P\{ X_n \in A |X_0=x\} \geq \lambda \varphi(A),
\end{eqnarray}
for all $x \in {\cal R}$ and $A \subset {\cal R}$.
Under (\ref{min}), Athreya and Ney \cite{Athreya1978}, and Nummelin \cite{Nummelin1978} show that $X_n$ admits a regenerative
scheme with i.i.d. interregeneration times for an augmented Markov chain, which is called the
``split chain''. It is known that under irreducibility and aperiodicity assumption, $w$-uniform ergodicity implies that $\{X_n,n \geq 0\}$ is Harris recurrent.

Denote $S_n:= \log LR_n(\theta),$ where $LR_n(\theta)$ is defined in (\ref{lr}).
Let $\varrho$ be the first time $(>0)$ to reach the atom of the split chain, and define
$u(\alpha,\zeta)= E_{\nu} e^{\alpha S_\varrho-\zeta \varrho}$ for $\zeta \in {\bf R}$, where $\nu$ is
an initial distribution on ${\cal X}$. Assume that
\begin{eqnarray}\label{mom}
W:= \{(\alpha,\zeta):u(\alpha,\zeta) < \infty\}~{\rm is~an~open~subset~on}~{\bf R}^2.
\end{eqnarray}
Denote $\zeta_1=\zeta_1(\theta):= \log LR_1(\theta).$
Ney and Nummelin \cite{Ney1987} shows that $D=\{\alpha:u(\alpha,\zeta) < \infty~{\rm for~some}~\zeta\}$ is an open set
and that for $\alpha \in D$, the transition kernel $\hat{P}_{\alpha}(x,A ) = E_{x}\{ e^{\alpha \zeta_1} I_{\{X_1 \in A \}}\}$ has a maximal simple real
eigenvalue $e^{\Psi(\alpha)}$, where $\Psi(\alpha)$ is the unique solution
of the equation $u(\alpha,\Psi(\alpha))=1$, with corresponding eigenfunction
$r^*(x;\alpha):= E_{x} \exp\{\alpha S_\varrho - \Psi(\alpha) \varrho\}$.
For a measurable subset $A \in {\cal B}({\cal X})$ and $x \in {\cal X}$, define
\begin{eqnarray}\label{eignmea}
{\cal L}(A ;\alpha) &=&  E_{\nu}\bigg[\sum_{n=0}^{\varrho -1} e^{\alpha S_n- n \Psi(\alpha)}
I_{\{X_n \in A \}}\bigg], \\
{\cal L}_{x}(A;\alpha) &=&  E_{x}\bigg[\sum_{n=0}^{\varrho -1} e^{\alpha S_n- n \Psi(\alpha)}
I_{\{X_n \in A \}}\bigg].
\end{eqnarray}

For each $\theta \in J$, denote $K({P}^{\theta}, {P}^{\theta_0})$ as the Kullback-Leibler
information numbers which will be defined precisely in (\ref{KL5.9}) of Section \ref{sec5}.

The following assumptions will be used throughout this paper.

\noindent
Condition C:

\noindent
C1. For each $\theta \in \Theta$, the Markov chain $\{X_n, n \geq 0 \}$
defined in (\ref{ssm}) and (\ref{ssmden}) is aperiodic, irreducible, and $V$-uniformly ergodic 
for some $V$ on
${\cal X}$, such that there exists $p \geq 1$,
\begin{eqnarray} \label{A12}
\sup_{x \in {\cal X} } E^{\theta}_{x}\bigg\{\frac{V(X_p)}{V(x)} \bigg\} < \infty.
\end{eqnarray}

\noindent
C2. For each $\theta \in \Theta$, assume $0< p_{\theta}(x,x')< \infty$ for all
$x,x' \in {\cal X}$, and $0< \sup_{x \in {\cal X}}
f(y;\theta|x)< \infty,$ for all $y \in {\bf R}^d$.
Denote $h_\theta(Y_1)= \sup_{x_{0} \in {\cal X}} \int p_{\theta}(x_{0},x_1)
f(Y_1;\theta|x_1) m(d x_1)$, and assume there exists $p \geq 1$ as in C1 such that
\begin{eqnarray}
&~& \sup_{x \in {\cal X}} E^{\theta}_{x} \bigg\{ \log \bigg( h_{\theta}(Y_1)^p
\frac{V(X_p)}{V(x)}  \bigg) \bigg\}  < 0 \label{3.1a}, \\
&~& \sup_{x \in {\cal X}} E^{\theta}_{x} \bigg\{ h_{\theta}(Y_1)
\frac{V(X_1)}{V(x)} \bigg\}  < \infty. \label{3.1b}
\end{eqnarray}

\noindent
C3. For each $\theta \in J$, assume $0< K({P}^{\theta},{P}^{\theta_0}) < \infty$. For each 
$\theta \in \Theta$, assume
\begin{eqnarray*}
&& \sup_{x_0 \in {\cal X} } | \int_{x_1 \in {\cal X}}
\int_{y \in {\bf R}^d} \pi_{\theta}(x_0)  p_\theta(x_0,x_1) f(y;\theta|x_1) \\
&&~~~~~~~ Q(dy) m(dx_1) |   < \infty.
\end{eqnarray*}

\noindent
C4. Assume (\ref{mom}) hold. Let $C$ be a measurable subset of ${\cal X} $ such that
\begin{eqnarray}\label{eigenpro}
{\cal L}(C;\alpha) < \infty~{\rm and}~{\cal L}_{x}(C;\alpha)< \infty~{\rm for~all~}
x \in {\cal X}.
\end{eqnarray}
Let $V: {\cal X} \to [1,\infty)$ be a measurable function such that for some $0 < \beta < 1$ and
$K > 0$,
\begin{eqnarray}
&& E_{x} [e^{\alpha \zeta_1- \Psi(\alpha)} V(X_1)] \leq (1 - \beta) V(x)~\forall~x \notin C, \label{eignmeapro1} \\
&& \sup_{x \in C} E_{x} [e^{\alpha \zeta_1- \Psi(\alpha)} V(X_1)] = K < \infty~\label{eignmeapro2} \\
&& {\rm and}~\int V(x)\varphi(dx) < \infty, \nonumber
\end{eqnarray}
where $\varphi$ is defined in (\ref{min}).\\

REMARK 2: C1 is an ergodic condition for the underlying Markov chain.
%considerable weak than the uniformly recurrent condition.
%A1 of Jensen and Petersen (1999), and that of Douc and Matias (2001).
The weighted mean contraction
property (\ref{3.1a}) and the finite weighted mean average property (\ref{3.1b}), appeared in C2,
guarantee that the induced Markovian iterated random functions system
satisfies uniformly ergodic condition with respect to a given norm. In Section \ref{secEA}, we show
that several interesting models satisfy these conditions.
%Note that both C1 and C2 are general enough to cover several interesting examples in Section 5.
C3 is a constraint of the Kullback-Leibler information numbers and a standard moment condition.
Note that positiveness of the Kullback-Leibler information numbers is not at all restrictive, since
it holds whenever the probability density functions of ${P}^{\theta}$ and ${P}^{\theta_0}$ do not
coincide almost surely. The finiteness condition is quite natural and holds in most cases.
C4 ensures the finiteness of the eigenfunction
$r(x;\alpha)$ and the eigenmeasure ${\cal L}(A ;\alpha)$, cf. Theorem 4 of Chan and Lai 
\cite{Chan2003}. These properties are useful for defining the exponential embedding in (\ref{et}) and (\ref{lrkf}) below. \\
%\cdf{uniform ergodic with respect to a given norm; Harris recurrent?}\\

The next theorem establishes second order approximation of the weighted SRP rule. 

%\cdf{investigate the detection delay at time 1 (not the worst case detection delay). Check Theorems 8.3.1 and 8.3.2, and 8.5.2.
%Sol. It is an equalizer rule in the sense that the conditional average delay to detection is a constant in the change point.
%For $C(\theta)$, check (8.389)--(8.391) for specific weight.}

\begin{theorem}\label{thm2}
Let $Y_1,\ldots,Y_n$ be a sequence of random
variables from a state space model $\{Y_n, n \geq 1\}$ satisfying conditions {\rm C1-C4}.
Suppose $F'(\theta)=dF(\theta)/d \theta$ exists, positive and continuous in an open neighborhood of $\theta \in \Theta$. Assume that ${S}_1$ is nonarithmetic with respect to
${P}_{\infty}^\theta$ and ${P}^\theta_1$. 
Then for given $x_0  \in {\cal X} $, as $b \rightarrow \infty$
\begin{eqnarray}\label{ae}
&~& {E}_1^\theta (N_b(F)|X_0=x_0) \\
&=& \frac{1}{K({P}^{\theta},{P}^{\theta_0})}
\bigg(b + \frac{1}{2}\log \frac{b}{K({P}^{\theta},{P}^{\theta_0})} + C(\theta) \bigg) + o(1),
\nonumber
\end{eqnarray}
where $C(\theta)$ will be defined precisely in {\rm (\ref{aec})} of Section \ref{sec5}.
\end{theorem}

The proof of Theorem \ref{thm2} is given in Section \ref{sec6}.

%This approximation is important for practical applications when the value of
%$K({P}^{\theta_1},{P}^{\theta_0})$ is moderate. Numerical computation of the constant $C(\theta)$ is an interesting and challenge task.
%For small value of $K({P}^{\theta_1},{P}^{\theta_0})$ the
%overshoot can be neglected, and the first order approximation in (\ref{ae}) will be reasonably accurate. 

The next theorem establishes asymptotic optimality of the weighted SRP rule.
%$0 < K({\mathbb P}^{\theta_0},{\mathbb P}^{\theta}) < \infty$,
%and ${\mathbb E}_1^\theta |{\mathbb S}_1|^2 < \infty$,

\begin{theorem}\label{thm1}
 Let $Y_1,\ldots,Y_n$ be a sequence of random variables from
a state space model $\{Y_n, n \geq 1\}$ satisfying conditions {\rm C1-C4}.
%$\theta_0 < a < b < \infty$, and
Assume $\theta_0 \in \Theta$, and suppose that there exists $J \subset \Theta$ with $F(J) > 0$.
Assume that for all $\theta \in J \subset \Theta$, ${S}_1$ is nonarithmetic with respect to
${P}_{\infty}^\theta$ and ${P}^\theta_1$. Then for any given change point detection rule $N \in
{\cal C} := \{E_{\infty}^{\theta_0} N \geq 1/B\}$, we have
\begin{eqnarray}\label{ao}
&& \inf_{N \in {\cal C}} \sup_{1 \leq \omega < \infty} \sup_{\theta \in J} 2 K({P}^{\theta}, 
{P}^{\theta_0}) E_{\omega}^{\theta} (N - \omega | N \geq \omega) \nonumber \\
&\geq&  2 b + \log b + O_\theta(1), 
\end{eqnarray}
where $ \limsup_{b \to \infty} \sup_{\theta \in J} |O_\theta(1)| < \infty$, and equality is attained by the weighted SRP rule.
\end{theorem}

The proof of Theorem \ref{thm1} is given in Section \ref{sec7}.

\def\theequation{3.\arabic{equation}}
\setcounter{equation}{0}
\section{Examples and applications}\label{secEA}

In this section, we demonstrate the application of our results to models of general state Markov models and
linear state space models, which are commonly used in practice for change point detection, cf. Tartakovsky et al. \cite{TNBbook14}. \\

\noindent
Example 1. General state Markov models

When $Y_n$ equals $X_n$ in (\ref{ssm}), one has a general state Markov chain.
Under the uniform recurrent condition for the Markov chain, and using the characterization of the
Kullback-Leibler distance $K(\theta,\theta_0):=\int_{x \in {\cal X}} \pi_{\theta}(x) \int_{x' \in {\cal X}}  p_{\theta}(x,x')  \log \frac{p_{\theta}(x,x')}{p_{\theta_0}(x,x')}dx'dx$,
Lai \cite{Lai1998} investigates the optimality property of generalized CUSUM rule under error probability constraint.
In this paper, we prove that the SRP rule is second-order asymptotic optimal, and present 
an asymptotic expansion of the average run length under conditions C1-C4. Note that the $V$-uniformly 
ergodic condition appeared in
C1 is weaker than the uniform recurrent condition, and covers several interesting examples. For instance an $AR(1)$ model with normal
innovation is $V$-uniformly ergodic with $V(x)=|x| +1$ [cf. pages 380 and 383 of Meyn and Tweedie \cite{Meyn2009}];
while it does not satisfy the assumption of the transition density function
$p_{\theta_1}(\cdot,\cdot)$ is uniform recurrent, in the sense that there exist
$c_2>c_1>0,~m \geq 1$ and a probability measure $\mu^*$ on ${\cal X}$ such that
$\displaystyle c_1 \mu^*(A) \leq P \{X_m \in A| X_0=x\} \leq c_2 \mu^*(A)$
for all measurable subsets A and all $x \in {\cal X}$.

To discuss condition C, appeared in Section \ref{secAO}, for this model.
Suppose that $\{X_n, n \geq 0\}$ is a Markov chain with transition density function
$p_{\theta_0}(\cdot, \cdot)$ for $n < \omega$ and $p_{\theta}(\cdot, \cdot)$ for
$n \geq \omega$, with respect to some $\sigma$-finite measure $m$ on the state space ${\cal X}$.
%In this case $g(Y_{k-1},Y_k)$ reduces to $\log \{p_{\theta}(X_{k-1}, X_{k})/ p_{\theta_0}(X_{k-1},X_{k})\}$.
Condition C1 requires that $\{X_n, n \geq 0\}$ is $V$-uniformly ergodic. 
By choosing $p=1$, $(\ref{A12})$
reduces to $\sup_{x}E_x^{\theta}V(X_1)/V(x) < \infty$. Condition C3 reduces to that for each $\theta \in \Theta$,
$0< p_\theta(x,y)  < \infty,$ for all $x,y \in {\cal X}$, which is also required
in C2. Note that $h(Y_1)$, used in (\ref{3.1a}) and (\ref{3.1b}),
reduces to $\sup_{x_0} \int p_\theta(x_0,x_1) \pi_\theta(x_1) m(dx_1)$.
Condition C4 reduces to a condition involves $X_0$ and $X_1$ only; see conditions (W1) and (W2) in
Chan and Lai \cite{Chan2003}.

One can show that many practical used models satisfy condition C. 
For instance, we consider an $AR(1)$ model
$X_n = \alpha X_{n-1} + \varepsilon_n,$ where $|\alpha| < 1$,  and $\varepsilon_n$ are independent and identically
distributed standard normal random variables. Under the normal errors assumption,
it is straightforward to check that C1 and C3 hold. To check condition C2, we only show that (\ref{3.1a}) holds
since the verification of (\ref{3.1b}) is the same. Note that $X_1$ has stationary distribution $N(0,a^2)$
with $a = 1/(1-\alpha^2)$. Observe that $Y_1=X_1$ and $h_\theta (Y_1)$ reduces to
\begin{eqnarray*}
&~& \sup_{x \in {\bf R}}   \int_{- \infty}^{\infty}
 \frac{\exp\{ -(y - \alpha x)^2/2 \}}{\sqrt{2 \pi}}  \frac{\exp \{ -y^2/2a^2 \} }{\sqrt{2 \pi}} dy  \\
&=&  \sup_{x \in {\bf R}} \frac{1}{\sqrt{2 \pi(1+a^2)}} \exp \bigg\{ -\frac{\alpha^2 x^2}{2(1+a^2)} \bigg\}  \int_{- \infty}^{\infty} \frac{\sqrt{1+a^2}}{\sqrt{2 \pi}a}  \\
&&~~~\times \exp\bigg\{ - \frac{1+a^2}{2 a^2}
\bigg(y - \frac{a^2\alpha x}{1+a^2}\bigg)^2 \bigg\}  dy \\
&=& \frac{1}{\sqrt{2 \pi(1+a^2)}} \sup_{x \in {\bf R}}  \exp \bigg\{ -\frac{\alpha^2 x^2}{2(1+a^2)} \bigg\}
=\frac{1}{\sqrt{2 \pi(1+a^2)}}.
\end{eqnarray*}
Consider $p=1$, a simple calculation leads that
\begin{eqnarray}\label{rem2.2a}
&& \sup_{x_0 \in {\bf R} } E^\theta_{x_0} \bigg\{ \log \bigg( h_\theta(X_1)
\frac{V(X_1)}{V(x_0)}  \bigg) \bigg\} \\
& <& \log  \sup_{x_0 \in {\bf R}} E^\theta_{x_0}\bigg\{ \frac{|\alpha x_0 + \varepsilon_{1}|+1}
{\sqrt{2 \pi(1+a^2)} (|x_0|+1)} \bigg\} \nonumber  \\
& \leq & \log  \sup_{x_0 \in {\bf R}} \bigg\{ \frac{|\alpha x_0| + E^\theta_{x_0} | \varepsilon_1|+1}
 {\sqrt{2 \pi(1+a^2)} (|x_0|+1)} \bigg\} \nonumber \\
& =& \log   \sup_{x_0 \in {\bf R}}\bigg\{ \frac{|\alpha x_0| + \frac{2}{\sqrt{2\pi}}+1}
{\sqrt{2 \pi(1+a^2)} (|x_0|+1)} \bigg\} < 0. \nonumber
\end{eqnarray}
This implies (\ref{3.1a}) hold. The verification of C4 is similar to Example 2 in Chan and Lai 
\cite{Chan2003}.

Next we consider the following example which involves change in the mean value
$\theta$ of a stable autoregressive sequence:
\begin{eqnarray}\label{arp}
 X_n = \sum_{k=1}^p a_k X_{n-k} + v_k + (1- \sum_{k=1}^p a_k) \theta,
\end{eqnarray}
where $a_1,\ldots,a_p$ are autoregressive coefficients and $v_k$ is a
Gaussian sequence with zero mean and variance $\sigma^2$. By Theorem
16.5.1 of Meyn and Tweedie \cite{Meyn2009}, $X_n$ defined in $(\ref{arp})$ is a
$V$-uniformly ergodic Markov chain with $V(x)=x^2+1$. It is easy to see condition C1 holds.
Since the verification of C2 can be done as that in (\ref{rem2.2a}), we will not repeat it here. Note that the assumption of normal distributed innovation (with mean zero and finite
variance $\sigma^2$) implies that the moment
condition C3 holds. The verification of condition C4 is similar to Example 2 in Chan and Lai 
\cite{Chan2003}. Note that this example can be generalized to the case of random coefficient
autoregression appeared on page 404 of Meyn and Tweedie \cite{Meyn2009}. \\

\noindent
Example 2. Linear state space models

Consider the stochastic system
\begin{eqnarray}
X_{n+1}&=& F X_n + G u_n + \delta_n, \label{lsma} \\
\|F\| &=& \sup\limits_{\|x\|=1}\|Fx\|<1, \nonumber \\
Y_{n}&=& H X_n + J u_n + \varepsilon_n, \label{lsmb} %\|H\|=\sup\limits_{\|x\|=1}\|Hx\|<1,
\end{eqnarray}
in which the unobservable state vector $X_n$, the input vector $u_n$, and the measurement vector
$Y_n$ have dimensions $p,q,$ and $r$, respectively, and $\delta_n,~\varepsilon_n$ are independent
Gaussian vectors with zero means and $cov(\delta_n)=\Sigma_1,~cov(\varepsilon_n)=\Sigma_2$. We assume
$G, J, \Sigma_1$ and $\Sigma_2$ are given, and the unknown parameter is $(F,H)^t,$ where $t$ denotes transpose.
The problem of additive change point detection can be found
in Tartakovsky et al. \cite{TNBbook14} and Lai \cite{Lai1998}. Here we consider the problem of 
nonadditive change.
Suppose at an unknown time $\omega$ the system undergoes some change in the sense that the parameter
is changed from $\theta_0$ to $\theta$, where $\theta_0$ is given while $\theta \in J \subset \Theta$ is unknown.
Here we consider $\theta$ is one dimensional unknown parameter, which can be one of the component in $(F,H)^t$,
the other parts are treated as nuisance parameters.

Let $\hat{H}_n$ and $\hat{F}_n$ be the estimators of $H$ and $F$, respectively.
The Kalman filter provides a recursive algorithm to compute the conditional expectation
$\hat{X}_{n|n-1}$ of the state $X_n$ given the past observations $Y_{n-1}, u_{n-1}, Y_{n-2}, u_{n-2},\ldots.$
The innovations $e_{n} = Y_n - \hat{H}_n \hat{X}_{n|n-1}- \hat{J}_n u_n$ are independent zero-mean Gaussian
vectors with $cov(e_n)=V_n$ given recursively by
\begin{eqnarray}\label{vn5}
V_{n} = \hat{H}_n P_{n|n-1}\hat{H}^t_n+\Sigma_2,
\end{eqnarray}
where
\begin{eqnarray}\label{pn}
 &~& P_{n+1|n} \\
 &=& \hat{F}_n (P_{n|n-1}-P_{n|n-1}\hat{H}^t_n
V_n^{-1} \hat{H}_n P_{n|n-1}) \hat{F}^t_n + \Sigma_1. \nonumber
\end{eqnarray}
When the parameter $\theta=\theta_0$, the innovations $e_n^0$ are independent Gaussian
vectors with covariance matrices $V_n^0$, and means $\mu_n^0=E(e_n^0)$ for $n \leq \omega$,
while when the parameter is changed to $\theta \in J$,
the innovations $e_n^\theta$ are independent Gaussian vectors with covariance matrices $V_n^\theta$,
and means $\mu_n^\theta=E(e_n^\theta)$ for $n \geq \omega$. Consider the weighted likelihood
\begin{eqnarray}\label{mlik}
LR_n^k(F) = \displaystyle \int_{\theta \in J} \prod_{l=k}^n
 \frac{f(e_l^\theta/\sqrt{V_l^\theta})}{f(e_l^0/\sqrt{V_l^0})}dF(\theta),
\end{eqnarray}
where $f(s)=e^{-||s||^2/2}/(2\pi)^{d/2}$ denotes the $d$-dimensional standard normal density, $d=p+r$,
and assume that the matrix whose inverse appears in $(\ref{mlik})$ is nonsingular.

To illustrate the computation of (\ref{mlik}), we consider a simple case that there is only a
one-dimensional unknown parameter $H=\theta \in J = (0,1)$ and $F(\theta)$ is uniform distributed on
$(0,1)$. Let $\theta + a_n=\mu_n$ for given $a_n$, and denote $\sigma_{n,\theta}^2$ as $V_n^\theta$.
When $\theta=\theta_0$, simply denote $a_l=0$ and $\sigma_{n,0}^2$ as $\sigma_{n,{\theta_0}}^2$.
That is $e_l^\theta \sim N(\theta+a_l,\sigma_{l,\theta}^2)$, and $e_l^0 \sim N(0,\sigma_{l,0}^2)$.
Then a simple calculation leads that
\begin{eqnarray*}
&~& {\rm (\ref{mlik})} \\
&=& \displaystyle \int_0^1 \prod_{l=k}^n
 \bigg[\exp \bigg\{-\frac{(e_l^{\theta}-(\theta+a_l))^2}{2 \sigma_{l,\theta}^2 }
 +\frac{(e_l^0-\theta_0)^2}{2\sigma_{l,0}^2}\bigg\} \bigg]d\theta \\
&=&\displaystyle\exp\bigg\{\sum_{l=k}^n \bigg[\frac{(e_l^0-\theta_0)^2}{2\sigma_{l,0}^2}-
\frac{(e_l^{\theta}-a_l)^2}{2\sigma_{l,\theta}^2}\bigg] \\
&~& +\frac{\displaystyle\bigg(\sum_{l=k}^n
(e_l^{\theta}-a_l)\bigg)^2}{2\alpha}\bigg \}\cdot\frac{\sqrt{2\pi}}{\sqrt{\alpha}}\bigg(\Phi(b)-\Phi(a)\bigg),
\end{eqnarray*}
where $\alpha= \sum_{l=k}^n 1/\sigma_{l,\theta}^2,~a=-\displaystyle\sum_{l=k}^n
\frac{e_l^{\theta}-a_l}{\sqrt{\alpha}}$, $b=\sqrt{\alpha}-\displaystyle\sum_{l=k}^n
\frac{e_l^{\theta}-a_l}{\sqrt{\alpha}}$, and $\Phi(\cdot)$ is the cumulative distribution
function of standard normal random variable.

Without assuming prior knowledge for the parameter after change $\theta$ and the change time $\omega$,
the weighted SRP change point detection rule, defined in Section \ref{secAO}, has the form
\begin{eqnarray}\label{ssmarp}
N_b &=& \inf \left \{ n: \sum_{k=0}^n LR_n^k(F)  \geq B \right \} \\
&=& \inf \left \{ n: \log \sum_{k=0}^n LR_n^k(F)  \geq b \right \}, \nonumber
\end{eqnarray}
where $B> 0$ is a given threshold and $b=\log B$.

To check the regularity condition C hold, we assume that there is no input vector $u_n$ for simplicity.
We first consider condition C1.
Note that $Y_n$ are independent for given $X_n$, therefore the weight function $V$ depends on $X_0$ only and
one can choose $V(x)=e^{\gamma\|x\|}$ for some $\gamma$ to be specified later.
Let $C=\{\mu:\|\mu\|\leq N\}$, and denote $\lambda$ as the Lebesgue measure on ${\bf R}^p$.
Recall that $\delta_1$ has normal density function $\phi$ with zero mean vector
and variance-covariance matrix $cov(\delta_1)=\Sigma_1$, which is positive and continuous,
and this implies $\eta:=\inf\{\phi(\delta-Fx):x \in C$ and $Fx + \delta \in C\}>0$.
Since $P\{F x_1 + \delta_1 \in d\delta \}\geq \phi(\delta-Fx) d\delta$,
we have for all $x \in {\bf R}^p$,
$P_{x} \{ X_1 \in A  \}\geq \delta I_{\{x \in C\}} \lambda(A \cap C)$, and therefore the minorization condition holds with $h(x)=\delta \lambda(C) \times I_{\{x \in C\}}$.
Under the normal error assumptions, it is easy to see that (\ref{A12}) and C3 hold.

To check condition C4 hold. Let
\begin{eqnarray}\label{5lr}
& \zeta_1 := \zeta_1 (\theta) \\
:=& \log \frac{\int_{x_0,x_1 \in {\cal X}} \pi_{\theta}(x_0)
    p_{\theta}(x_{0},x_1) f(Y_1;\theta|x_1)m(dx_1)m(dx_0)}
{\int_{x_0,x_1 \in {\cal X}}  \pi_{\theta_0}(x_0)
    p_{\theta_0}(x_0,x_1) f(Y_1;\theta_0|x_1)m(dx_1) m(dx_0)}, \nonumber
\end{eqnarray}
where $\pi_{\theta}(x_0)$ is the $p$-variate normal density function with zero mean vector and variance-covariance
matrix $\Sigma_1/(1-||F||)$, $p_{\theta}(x_{0},x_1)$ is the $p$-variate normal density function with mean
vector $F x_0$ and variance-covariance matrix $\Sigma_1$, and $f(Y_1;\theta|x_1)$ is the $p$-variate normal density function with mean vector $H x_1$ and variance-covariance matrix $\Sigma_2$.
Denote the conditional distribution of $\zeta_1$ given $(X_1,Y_1)$ has the form
$F_{(X_1,Y_1)}$. Since $\zeta_n = g(Y_n)$ for some $g$ by (\ref{5lr}), $F_{(X_1,Y_1)}$
degenerates to $F_{Y_1}$. By (\ref{lsma}), (\ref{lsmb}) and (\ref{5lr}), it is easy to see that for any given $\alpha \in {\bf R}$, there exists a positive constant $\rho_\alpha$ such that
\begin{eqnarray}\label{527}
\int e^{\alpha g(s_1)}dF_{s_1}(s_1) \leq \exp\{\rho_\alpha\|s_1\|\}~
{\rm for~all~} s_1 \in {\bf R}^r.
\end{eqnarray}
%Furthermore we have that for every compact subset $C$ of
%${\bf R}$, there exists a finite measure $\nu_C$ with compact support $K_C$ such that
%\begin{eqnarray}\label{528}
% \inf\limits_{x,y \in C}F_{x,y}(B) \geq \nu_C(B)~~~~~~~{\rm for~all}~B\subset K_C.
%\end{eqnarray} Moreover,by (\ref{527}),
This implies that
\begin{eqnarray}
&~& E_{x} [e^{\theta \zeta_1}V(X_1)] \\
&\leq& E \exp\{\rho_\theta(\|Hx +  \varepsilon_1 \|)+\gamma\|Fx + \delta_1\|\} \nonumber \\
&\leq& \Lambda(\rho_\theta+\gamma)\exp\{\big(\rho_\theta(1+\|H\|)+\gamma\|F\|\big)\|x\|\}. \nonumber
\end{eqnarray}
Since $\|F\|<1$, we can choose $\gamma$ large enough so that $2\rho_\theta+\gamma\|H\|<\gamma$, and then
(\ref{eignmeapro1}) is satisfied if $N$ is chosen large enough.
Since $C$ is compact and $\lambda(\cdot\cap C)$ has support $C$, (\ref{eignmeapro2}) also holds for sufficiently
large $L$.

Finally we need to verify C2 hold. For simplicity, let $p=r=d$ in
(\ref{lsma}) and (\ref{lsmb}). After normalization, we may assume the variance parts in $\Sigma_1$ and $\Sigma_2$ are both equal
to $1$. That is, define
\[
\Sigma_1 = \left(
\begin{array}{cccc}
1        & \rho_1^1 & \cdots & \rho_d^1 \\
\rho_1^1 &        1 & \cdots & \rho_{d-1}^1 \\
\vdots   & \vdots   & \ddots & \vdots \\
\rho_d^1 & \rho_{d-1}^1 & \cdots & 1
\end{array} \right)\]
and
\[
\Sigma_2 = \left(
\begin{array}{cccc}
1        & \rho_1^2 & \cdots & \rho_d^2 \\
\rho_1^2 &        1 & \cdots & \rho_{d-1}^2 \\
\vdots   & \vdots   & \ddots & \vdots \\
\rho_d^2 & \rho_{d-1}^2 & \cdots & 1
\end{array} \right). \]

Let
\[
x_0=x = \left( \begin{array}{c}
x_1\\
x_2\\
\vdots\\
x_d \end{array} \right), ~
x_1=x' = \left( \begin{array}{c}
x'_1\\
x'_2\\
\vdots\\
x'_d \end{array} \right),~
y = \left( \begin{array}{c}
y_1\\
y_2\\
\vdots\\
y_d \end{array} \right),
\]
\[
z=Hx'=
\left(
\begin{array}{cccc}
h_{11}     & \cdots & \cdots & h_{1d} \\
\vdots &   \cdots      & \cdots & \cdots \\
\vdots   & \vdots   & \ddots & \vdots \\
h_{d1} & \cdots & \cdots & h_{dd}
\end{array} \right)
\left(
\begin{array}{c}
x'_1\\
x'_2\\
\vdots\\
x'_d \end{array} \right),\]

\[\mu = FX=
\left(
\begin{array}{cccc}
\alpha_{11}     & \cdots & \cdots & \alpha_{1d} \\
\vdots &   \cdots      & \cdots & \cdots \\
\vdots   & \vdots   & \ddots & \vdots \\
\alpha_{d1} & \cdots & \cdots & \alpha_{dd}
\end{array} \right)
\left( \begin{array}{c}
x_1\\
x_2\\
\vdots\\
x_d \end{array} \right),
\]
\[ \mu^* = (\Sigma_1^{-1}+H^t\Sigma_2^{-1}H)^{-1}(\Sigma_1^{-1}\mu+H^t\Sigma_2^{-1}y) \]
and
\[\Sigma^{*-1}=\Sigma_1^{-1}+H^t\Sigma_2^{-1}H. \]
Denote $|\Sigma|$ as the determinant of the matrix $\Sigma$. Then a simple calculation leads that
\begin{eqnarray*}
&~& \displaystyle   \int_{-\infty}^{\infty} \cdots \int_{-\infty}^{\infty}
\frac{\exp\{-\frac{1}{2} (x'-\mu)^t\Sigma_1^{-1}(x'-\mu)\}}{(2\pi)^{d/2}|\Sigma_1|^{1/2}} \\
&~& \times  \frac{exp\{-\frac{1}{2}  (y-z)^t\Sigma_2^{-1}(y-z)\}}{(2\pi)^{d/2}|\Sigma_2|^{1/2}} dx'_1 \cdots dx'_d \\
&=& \displaystyle\frac{|\Sigma^*|^{1/2}}{(2\pi)^{d/2}|\Sigma_1|^{1/2}|\Sigma_2|^{1/2}}
\exp \bigg\{ \mu^t\Sigma_1^{-1}\mu+ s^t\Sigma_2^{-1}y  \\
&~& - \bigg[(\Sigma_1^{-1}+H^t\Sigma_2^{-1}H)^{-1}
 (\Sigma_1^{-1}\mu+H^t\Sigma_2^{-1}y)\bigg]^t\\
&~& \displaystyle \times ~\Sigma^{*-1} \bigg[(\Sigma_1^{-1}+H^t\Sigma_2^{-1}H)^{-1}(\Sigma_1^{-1}\mu+H^t\Sigma_2^{-1}y)\bigg] \bigg\}.
\end{eqnarray*}
Note that
\begin{eqnarray*}
&& \Sigma^{*-1}=\Sigma_1^{-1}+H^t\Sigma_2^{-1}H \Longrightarrow \Sigma^* = (\Sigma_1^{-1}+H^t\Sigma_2^{-1}H)^{-1} \\
&&\Longrightarrow |\Sigma^*|=|(\Sigma_1^{-1}+H^t\Sigma_2^{-1}H)^{-1}|. 
\end{eqnarray*}
%Taking $\sup_{x \in {\bf R}^d}$ to get
Therefore
\begin{eqnarray}\label{ex2a}
&~& h(y) \\
&=& \displaystyle \sup_{x \in {\bf R}^d}  \int_{-\infty}^{\infty} \cdots \int_{-\infty}^{\infty}
\frac{\exp\{-\frac{1}{2}(x'-\mu)^t\Sigma_1^{-1}(x'-\mu)\}}{(2\pi)^{d/2}|\Sigma_1|^{1/2}} \nonumber \\
&&~\times  \frac{\exp\{-\frac{1}{2} (y-z)^t\Sigma_2^{-1}(y-z)\}}{(2\pi)^{d/2}|\Sigma_2|^{1/2}} dx'_1 \cdots dx'_d \nonumber \\
&=& \displaystyle  \frac{1}{(2\pi)^{d/2}|\Sigma_1|^{1/2}|\Sigma_2|^{1/2}}\cdot
|(\Sigma_1^{-1}+H^t\Sigma_2^{-1}H)|^{-1/2}. \nonumber
\end{eqnarray}

Assume $\frac{|(\Sigma_1^{-1}+H^t\Sigma_2^{-1}H)|^{-1/2}}{(2\pi)^{d/2}|\Sigma_1|^{1/2}|\Sigma_2|^{1/2}}
 =a <1 $. % and let $b =(1-|\alpha|^{2p})/(1-|\alpha|^2)$.
A simple calculation leads that
\begin{eqnarray}\label{rem2.2}
&~&\sup_{x_0 \in {\bf R}^d }
E^\alpha_{x_0} \bigg\{ \log \bigg( h(Y_1)^p
\frac{w(X_p)}{w(x_0)}  \bigg) \bigg\}  \\
&=&  \sup_{x_0 \in {\bf R}^d} E^\alpha_{x_0} \log \bigg\{ \frac{a^p\exp\{\gamma(\alpha^p \|x_0\|
 + \sum_{k=0}^{p-1} \alpha^k \varepsilon_{p-k})\}} {\exp\{\gamma \|x_0\| \}} \bigg\} \nonumber \\
&=&   \sup_{x_0 \in {\bf R}^d} E^\alpha_{x_0} \bigg\{ \gamma \alpha^p \|x_0\|
+ \sum_{k=0}^{p-1} \alpha^k \varepsilon_{p-k} - \gamma \|x_0\| + p \log a \bigg\} \nonumber \\
&=& p \log a < 0. \nonumber
\end{eqnarray}
%Chooding $p>\log( 2/(1-|\alpha|^2) +1)/\log (2\pi \sqrt{4-(\rho_1+\rho_2)^2})$, we have (\ref{rem2.2})$<0$.
This implies (\ref{3.1a}) hold. By using the same argument, we have (\ref{3.1b}) hold.

To illustrate (\ref{ex2a}) and (\ref{rem2.2}), we consider a simple case of $d=2$. Denote
$x=\left(\begin{array}{c}
x_1\\
x_2\end{array}\right),~
x' =\left(\begin{array}{c}
x'_1\\
x'_2\end{array}\right),~
\mathbf{\mu} = \left(
\begin{array}{c}
\alpha_1 x_1 \\
\alpha_2 x_2\end{array}\right),~
y=\left(
\begin{array}{c}
y_1\\
y_2\end{array}\right),$
$\Sigma_1 = \left(
\begin{array}{cc}
 1&\rho_1\\
 \rho_1&1
\end{array} \right)$ and
$\Sigma_2 = \left(
\begin{array}{cc}
 1&\rho_2\\
 \rho_2&1
\end{array} \right).$
Simple calculation leads that
\begin{eqnarray*}
\Sigma^{*-1} &=&
\frac{2-\rho_1^2-\rho_2^2}{(1-\rho_1^2)(1-\rho_2^2)} \\
&&~ \times \left(
\begin{array}{cc}
 1& \frac{-\rho_1(1-\rho_2^2)-\rho_2(1-\rho_1^2)}{2-\rho_1^2-\rho_2^2} \\
 \frac{-\rho_1(1-\rho_2^2)-\rho_2(1-\rho_1^2)}{2-\rho_1^2-\rho_2^2} &1
\end{array} \right),
\end{eqnarray*}
and
$$
\mu^*=\left( \begin{array}{c}
\frac{(\rho_1-\rho_2)(\mu_2-y_2)+(\rho_1^2+\rho_1\rho_2-2)s_1+(\rho_2^2+\rho_1\rho_2-2)\mu_1}{(\rho_1+\rho_2)^2-4}\\
\frac{(\rho_1-\rho_2)(\mu_1-y_1)+(\rho_1^2+\rho_1\rho_2-2)s_2+(\rho_2^2+\rho_1\rho_2-2)\mu_2}{(\rho_1+\rho_2)^2-4}
\end{array} \right).$$
Then $h(y)=\frac{1}{2\pi \sqrt{4-(\rho_1+\rho_2)^2}}$, and the condition reduces to
$|\rho_1+\rho_2| < \frac{\sqrt{16\pi^2 -1}}{2\pi}  \approx 1.994$.  \\

\def\theequation{4.\arabic{equation}}
\setcounter{equation}{0}
\section{Likelihood representation and exponential embedding}\label{secLR}

In this section, we investigate the weighted Shiryayev-Roberts change point detection rule (\ref{mlrk})-(\ref{msr}). Due to the change point detection rule involves $LR_n^k (\theta)$ defined in (\ref{lrka}),
we study the likelihood ratio $LR_n$ appeared in (\ref{lr}) first.
%as (\ref{lrk}) and (\ref{mlrk}) can be analyzed in the same manner.
A major difficulty for analysing the likelihood ratio (\ref{lr}) is its integral form. To overcome this obstacle,
we represent (\ref{lr}) as the ratio of $L_1$-norms of a Markovian iterated random functions system.
%This device has been proposed by Fuh (2004) to study efficient likelihood estimation in state space models.
%Here, we carry out the same idea to have a representation of the likelihood ratio $LR_n$.
Specifically, let
\begin{eqnarray}\label{mf}
{\bf H} &=& \{h|h:{\cal X} \to {\bf R}^+ ~{\rm is~}m{\rm -measurable},~  \\
&&~~~\int h(x) m(dx) < \infty~{\rm and~}\sup_{x \in {\cal X}} h(x)< \infty \},\nonumber
\end{eqnarray}
and define the variation distance between any two elements $h_1,h_2$ in ${\bf H}$ by
\begin{eqnarray}
d(h_1,h_2) = \sup_{x \in \cal X} |h_1(x)-h_2(x)|.
\end{eqnarray}

For $j=1,\ldots,n$, define the random functions ${\bf P}_\theta(Y_j)$
on ${\cal X} \times {\bf H}$ as
\begin{eqnarray}
&~&  {\bf P}_\theta(Y_0)h(x) \label{rf1} \\
&=& \int_{x' \in {\cal X}} f(Y_0;\theta|x') h(x') m(dx')~{\rm a~constant},  \nonumber  \\
&~& {\bf P}_\theta(Y_j)h(x) \label{rf2} \\
&= & \int_{x' \in {\cal X}} p_{\theta}(x,x') f(Y_j;\theta|x') h(x') m(dx'), \nonumber
\end{eqnarray}
and denote the composition of two random functions as 
\begin{eqnarray}\label{crf}
&~& {\bf P}_\theta(Y_{j+1})\circ{\bf P}_\theta(Y_{j})h(x) \\
&=& \int_{x_j \in {\cal X}} p_\theta(x,x_j) f(Y_j;\theta|x_j)  \bigg(\int_{x_{j+1} \in {\cal X}}
p_\theta(x_j,x_{j+1}) \nonumber \\ 
&~&~~~\times  f(Y_{j+1};\theta|x_{j+1}) h(y) m(dx_{j+1})\bigg) m(dx_j). \nonumber
\end{eqnarray}

%\cdf{backward; normalized; locally compact; separable; Harris recurrent, etc..}

Furthermore, let  
\begin{eqnarray}
{\bf M} &=& \{ M: {\bf H} \rightarrow {\bf H}| M {\rm~ is~a~linear} \\
~~&~&~~~~~~~~{\rm ~and~bounded~operator}~ 
	P^{\pi}_{\theta^*}-a.s. \}, \nonumber
\end{eqnarray}
be equipped with the operator norm $\|\cdot\|$ with respect to the sup-norm, i.e.
\begin{equation}
	\|M\| = \sup_{h\in {\bf H}:\|h\|_{\infty}=1}\|M(h)\|_{\infty}.
\end{equation}
We define the iterated random functional system as
\begin{eqnarray}
	M_{\theta,0}(h) &=& {\bf P}_\theta(Y_0) h \\
	M_{\theta,n}(h) &=& F(Y_{n},M_{\theta,n-1})(h) \label{MIRFS1} \\
	&:=& \frac{M_{\theta,n-1}(P_{\theta}(Y_n)h)}{\int M_{\theta,n-1}(P_{\theta}(Y_n)1)(x) m(dx)},  
	\nonumber
% 	&= M_{-1}\circ P_{\theta,0}\circ P_{\theta,1} \circ \cdots \circ P_{\theta,n}(h)\label{MIRFS2}
\end{eqnarray}
for $n\geq 1.$ Note that 
\begin{eqnarray}
&~&	M_{\theta,n}(h)(x) \\
&=& \int_{x_n\in {\cal X}} h(x_n)p_{\theta}\left(X_0=x,X_n=x_n |
	Y_0,\cdots,Y_n\right)m{(dx_n)} \nonumber
\end{eqnarray} 

For $h \in {\bf M}$, let $\|h\|:= \int_{x \in {\cal X}} h(x) m(dx)$ be the $L^1$-norm
on ${\bf M}$ with respect to $m$.
Then, the likelihood ratio (\ref{lr}) can be represented as
\begin{eqnarray} \label{lrrf}
LR_n(\theta) = \frac{||{\bf P}_{\theta}(Y_n) \circ \cdots \circ {\bf P}_{\theta}(Y_1)
\pi_{\theta}||} {||{\bf P}_{\theta_0}(Y_n) \circ \cdots \circ {\bf P}_{\theta_0}(Y_1)
 \pi_{\theta_0}||}.
\end{eqnarray}
%Note that, for $j=1,\ldots,n,$ the integrate
%$p_\theta(x,y)f(\xi_j;\varphi_y(\theta)|\xi_{j-1})$ of ${\bf P}_\theta(
%\xi_j)$ in (2.6) and (2.8) represents $X_{j-1} =x$ and $X_j\in dy$,
%and $\xi_j$ is a Markov chain with transition probability density $f(\xi_j;\varphi_{y}(\theta)|\xi_{j-1})$ for
%given ${\bf X}$. By definition (2.1), $\{(X_n,\xi_n), n \geq 0 \}$ is a Markov chain, and
%this implies that ${\bf P}_\theta(\xi_j)$ is a sequence of Markovian iterated random
%functions system (see Section 3 for a formal definition). Therefore, by representation (2.8), $p_n(\xi_0,\xi_1,\ldots,\xi_n;\theta)$ is the $L_1$-norm of a Markovian iterated random functions system.

Let $\{(X_n,Y_n), n \geq 0 \}:=\{(X_n^\theta,Y^\theta_n), n \geq 0 \}$ be the Markov
chain defined in (\ref{ssm1}) and (\ref{ssmden}).
Abuse the notation a little bit, we denote $\theta=(\theta_0,\theta)$ because $\theta_0$ is given.
For each $n$, let
\begin{eqnarray} \label{rfs}
 M_n(\theta) &=& {\bf P}_\theta(Y_n) \circ \cdots \circ {\bf P}_\theta(Y_1)
 = (M_n(\theta_0), M_n(\theta)) \\
&=& \big({\bf P}_{\theta_0}(Y_n) \circ \cdots \circ {\bf P}_{\theta_0}(Y_1),
{\bf P}_{\theta}(Y_n) \circ \cdots \circ {\bf P}_{\theta}(Y_1) \big) \nonumber
\end{eqnarray}
be the Markovian iterated random functions system on ${\bf M}$ induced from (\ref{rf2}). Then $\{W_n^\theta, n \geq 0\}:=\{(X_n^\theta,M_n(\theta)), n\geq 0\}$ is a Markov chain
on the state space ${\cal X}  \times {\bf M}$, with transition probability kernel
\begin{eqnarray}\label{tk}
 {\mathbb P}^{\theta}((x_0,h), A  \times \Gamma) &:=& \int_{x_1 \in A} \int_{y \in B} I_\Gamma( {\bf P}_\theta(y) h)\\
&& \times p_\theta(x_0,x_1) f(y;\theta|x_1) Q(dy) m(dx_1) \nonumber
\end{eqnarray}
for all $x_0 \in {\cal X},~h \in{\bf M},~A \in {\cal B}({\cal X})$ and
$\Gamma \in {\cal B}({\bf M})$, where $I_{\Gamma}$ denotes the indicator function on the set $\Gamma$.
%The $n$-step transition kernel is denoted ${\mathbb P}^{n}$.
For $(x,h)\in {\cal X} \times {\bf M}$, let ${\mathbb P}_{(x,h)}$ be the probability
measure on the underlying measurable space under which $X_0= x,M_{0}=h$. 
The associated expectation is denoted ${\mathbb E}_{(x,h)}$,
as usual. For an arbitrary distribution $\nu$ on ${\cal X}  \times {\bf M}$,
we put ${\mathbb P}_{\nu}(\cdot) := \int {\mathbb P}_{(x,h)}(\cdot)\,
\nu(dx \times dh)$ with associated expectation ${\mathbb E}_{\nu}$.
We use ${\mathbb P}$ and ${\mathbb E}$ for probabilities and expectations,
respectively, that do not depend on the initial distribution.
Since the Markov chain $\{X_n, n \geq 0\}$ has transition
probability density and the iterated random function $M_1(\theta)$, defined in (\ref{rfs}), is driven by
$\{(X_n,Y_n), n \geq 0\}$, the induced transition probability
$\mathbb{P}(\cdot,\cdot)$ has a density with respect to $m  \times Q$.
Denote it as $\mathbb{P}$ for simplicity. According to Theorem 1(iii) in Fuh \cite{Fuh2006}, the stationary distribution of $\{W^{\theta}_n, n\geq 0\}$ exists, and denote it by $\Pi_{\theta}$.

Now the log-likelihood ratio can be written as an additive functional of the Markov chain
$\{W^\theta_n, n \geq 0\}$. That is
\begin{eqnarray}\label{lrn}
\log LR_n(\theta)=\sum_{k=1}^{n}g(W^\theta_{k-1},W^\theta_k),
\end{eqnarray}
where
\begin{eqnarray}\label{gk}
g(W^\theta_{k-1},W^\theta_k) 
&:=& \log \frac{||{\bf P}_{\theta}(Y_k) \circ \cdots \circ {\bf P}_{\theta}(Y_1)  \pi_\theta||}
{||{\bf P}_{\theta_0}(Y_k) \circ \cdots \circ {\bf P}_{\theta_0}(Y_1) \pi_{\theta_0}||}  \\
&& - \log \frac{||{\bf P}_{\theta}(Y_{k-1}) \circ \cdots \circ {\bf P}_{\theta}(Y_1) \pi_\theta||}
{||{\bf P}_{\theta_0}(Y_{k-1}) \circ \cdots \circ {\bf P}_{\theta_0}(Y_1) \pi_{\theta_0}||}. \nonumber
\end{eqnarray}

To analyze the weighted SRP change point detection rule in state space models,
we need first to construct an exponential embedding of the transition probability
operator for the induced Markov chain $\{W_n, n \geq 0\}$ with state space ${\cal W}:={\cal X} \times {\bf M}$,
%Here, $Y_n = Y_n^{\theta_1}$ and $\theta_1=(\theta_0,\theta_1)$ for given $\theta_1 \in J \subset \Theta$,
and then to represent the weighted likelihood ratio as an additive functional of the Markov chain $\{W_n, n \geq 0\}$.
To this end, we show that the induced Markov chain
$\{W_n, n \geq 0\}$ satisfies some required recurrent and ergodic conditions.

For any given two transition probability kernels $Q(w,A),K(w,A)$,
$w \in {\cal W}$, $A \in {\cal B}({\cal W})$, the $\sigma$-algebra of ${\cal W}$, and for all measurable functions
$h(w), w \in {\cal W}$, define $Q h$ and $Q K$ by
$Q h(w) = \int Q(w,dw') h(w')$ and $Q K(w,A)=\int K(w,dw')Q(w',A)$, respectively.
Let ${\cal N}$ be the Banach space of measurable
functions $h:{\cal W} \rightarrow {\bf C}$ (:= the set of complex numbers) with
norm $\|h\| < \infty$. We also introduce the Banach space ${\cal B}$
of transition probability kernels $Q$ such that the operator norm
$||Q||= \sup\{||Qg||;||g|| \leq 1\}$ is finite.

Denote by $P^n(y,A)=P\{W_n \in A|W_0=y\}$, the transition probabilities over $n$ steps. The kernel $P^n$ is a $n$-fold power of $P$. Define also the C\'{e}saro averages $P^{(n)} = \sum_{j=0}^n P^j/n$, where
$P^0=P^{(0)}=I$ and $I$ is the identity operator on ${\cal B}$.

\begin{definition}
A Markov chain $\{X_n,n \geq 0\}$ is said to be uniformly ergodic with
respect to a given norm $||\cdot||$, if there exists a stochastic kernel
$\Pi$ such that $P^{(n)} \rightarrow \Pi$ as $n \to \infty$ in the
induced operator norm in ${\cal B}$.
\end{definition}

\begin{definition}\label{def3}
Let $\omega:{\cal X} \to [1,\infty)$ be the weight function defined in {\rm C1}, and ${\bf M}$ be
defined in $(\ref{mf})$. For any measurable function $g: {\cal X} \times {\bf M} \rightarrow [1,\infty)$, define
$||g||_V := \sup_{(x,h)  \in {\cal X} \times {\bf M}}
  \frac{|g(x,h)|}{V(x)},$
 and
$ \|g\|_h := \sup_{x \in {\cal X}, h_1,h_2: 0 < d(h_1,h_2) \leq 1}
\frac{|g(x,h_1) - g(x,h_2)|}{(V(x) d(h_1,h_2))^{\delta}},$
for $0 < \delta <1$.
We define ${\cal H}$ as the set of $g$ on ${\cal X} \times {\bf M}$ for which $\|g\|_{V h} := \|g\|_V + \|g\|_h$ 
is finite, where $Vh$ represents a combination of the weighted variation norm and the bounded 
weighted H\"{o}lder's norm.
\end{definition}

\begin{theorem}\label{thm3}
Let $\{(X_n^\theta,Y^\theta_n), n \geq 0\}$ be the state space model given in {\rm (\ref{ssm})},
satisfying {\rm C1-C3}, where $\theta=(\theta_0, \theta) \in \Theta \times J$ is the unknown parameter.
Then the induced Markov chain $\{W_n^\theta,n \geq 0\}$ is
an aperiodic, irreducible and Harris recurrent Markov chain. Moreover, it is
uniformly ergodic with respect to the norm defined in {\rm Definition \ref{def3}}. Furthermore
there exist $a,C > 0,$ such that ${\mathbb E}_{w}(\exp\{a g(W_0,W_1)\}) \leq C < \infty$ for all $w \in {\cal W}.$
\end{theorem}

Since the proof is the same as those in Lemmas 3 and 4 of Fuh \cite{Fuh2006}, it is omitted. 

Next we define Laplace transform of
the transition operator and introduce the twisting probability 
measure for $\{W_n, n \geq 0\}$.
Denote $w:=(x,h)$ and $\tilde{w} := (x_0, \pi)$, where $x_0$
is the initial state of $X_0$ taken from $\pi(X_0)$.
Recall $g(W_0,W_1)$ defined in (\ref{gk}).
For given %$L_1$-bounded measurable functions $h$,
$w \in {\cal W}$, $A \times \Gamma \in \mathcal{B(W)}$, and $\alpha \in {\bf  R}$, define the linear operator $\hat{{\bf P}}_{\alpha}$ by
\begin{eqnarray}\label{ee}
 \hat{{\bf P}}_{\alpha} (w, A \times \Gamma) = {\mathbb E}_w \bigg\{e^{\alpha g(W_0, W_1)}
     I_{\{W_1 \in A \times \Gamma\}} \bigg\}.
%{\bf P} h(y) = {\mathbb E}_y \{h(Y_0,Y_0)\}.
\end{eqnarray}
Under conditions C1-C3, Theorem \ref{thm3} leads that 
$\{W_n, n \geq 0\}$ is
an aperiodic, irreducible and Harris recurrent Markov chain, and conditions in
Theorem 4.1 of Ney and Nummelin \cite{Ney1987} hold. Therefore, 
$\hat{{\bf P}}_{\alpha}$ has
a maximal simple real eigenvalue $\lambda(\alpha)$ with associated right eigenfunction $r(\cdot;\alpha)$ such that $\Lambda(\alpha) = \log \lambda(\alpha)$ is
analytic and strictly convex on 
${\cal D}=\{\alpha: \Lambda(\alpha) < \infty\}$.

Let $\tau$ be the first time $(>0)$ to reach the atom of the split chain for $\{W_n, n \geq 0\}$.
For each $w \in \mathcal{W}$ and $A \times \Gamma \in \mathcal{B(W)}$, define the left eigenmeasures
\begin{eqnarray}\label{a1}
 \ell(A \times \Gamma;\alpha)&:= {\mathbb E}_{\nu}
\bigg\{\sum\limits^{\tau -1}_{n=0}e^{\alpha S_n - n \Lambda(\alpha)} I_{\{W_n \in A \times \Gamma\}}\bigg\}, \\
 \ell_w (A \times \Gamma;\alpha)&:= {\mathbb E}_{w} \bigg\{\sum\limits^{\tau -1 }_{n=0}
e^{\alpha S_n - n \Lambda(\alpha)} I_{\{W_n \in A \times \Gamma \}}\bigg\}. \nonumber
\end{eqnarray}
Recall $S_n = \sum_{k=1}^n g(W_{k-1},W_k)$, and $g(W_{k-1},W_k)$ defined in (\ref{gk}) is an additive functional of the Markov chain $\{(W_{n-1},W_n), n \geq 1\}$.
Since $r(w;\alpha)^{-1}  \pi_\alpha(dw) = L_\alpha \ell(dw;\alpha)$ for some constant $L_\alpha$ [cf.
Ney and Nummelin  \cite{Ney1987}, page 581], the finiteness of $\ell(A \times \Gamma;\alpha)$
implies that $r(w;\alpha)>0$ uniformly for $w \in {\cal W}$. On
the other hand,
Theorem 4 of Chan and Lai \cite{Chan2003} establishes the finiteness of $\ell(A \times \Gamma;\alpha)$ and $\ell_w (A \times \Gamma;\alpha)$.

Denote $\theta:=(\theta_0, \theta) \in \Theta \times J$ as the parameter.
Assume $\theta= \Lambda'(\alpha)$ is a one to one function, so one can indifferently consider
$\theta$ to be a function of $\alpha$ or $\alpha$ a function of $\theta$. Here $'$ denotes derivative.
For simplicity, we replace
$\alpha$ by $\theta$ in (\ref{ee}), and let ${\cal D}=J$ here and in the sequel.
Then under conditions C1-C4, by using Theorem 1 of Ney and Nummelin \cite{Ney1987} and Theorem 4 of 
Chan and Lai \cite{Chan2003}, we have
$r(\cdot;\theta)$ is uniformly positive, bounded and analytic on $J$ for each $w \in {\cal W}$.
For $\theta \in J$, define the twisting transformation for
the transition probability of $\{W_n, n \geq 0\}$ as
\begin{eqnarray}\label{et}
{\mathbb P}^{\theta}(w,dw')= \frac{r(w';\theta)}{r(w;\theta)}
     e^{-\Lambda(\theta) + \theta g(W_0, W_1)} {\mathbb P}(w,dw').
\end{eqnarray}

For given $\theta \in J \subset \Theta \subset {\bf R}$,
let $\{W_n^{\theta}, n \geq 0\}$ be the Markov 
chain with transition kernel
${\mathbb P}^{\theta}$ and invariant probability $\Pi^{\theta}$. If the function $\Lambda(\theta)$
is normalized so that $\Lambda(0)=\Lambda^{'}(0)=0$, then 
${\mathbb P}={\mathbb P}^0$
is the transition probability of the Markov chain $\{W_n, n \geq 0\}$, with invariant probability $\Pi= \Pi^{0}$.

By making use of (\ref{et}) and repeat the same idea as (\ref{lrn}), we have representations for
\begin{eqnarray}\label{lrk}
 LR^k_n(\theta) = \exp\bigg(\sum_{l=k+1}^{n}g(W^\theta_{l-1},W^\theta_l)\bigg),
\end{eqnarray}
and
\begin{eqnarray}\label{lrkf}
 LR_n^k(F) %= \int_{\theta \in J}  LR^k_n(\theta) dF(\theta)
&=& \int_{\theta \in J}\frac{r(W_n;\theta)}{r(W_k;\theta)}
     \exp\bigg\{-(n-k) \Lambda(\theta) \nonumber \\
&& + \theta \sum_{l=k+1}^n g(W_{l-1}, W_{l})\bigg\} dF(\theta). 
\end{eqnarray}~

\def\theequation{5.\arabic{equation}}
\setcounter{equation}{0}
\section{Second order approximation of the weighted SRP detection rule}\label{sec5}

By using the same idea as that in Pollak \cite{Pollak1985} and Fuh \cite{Fuh2004b}, we
introduce a randomization on the initial $LR_n^0(\theta)$ for the Shiryayev-Roberts scheme, and call
it the Shiryayev-Roberts-Pollak (SRP) change point detection rule in state space models. Before that,
we need the following notations first.
%we consider a randomization for a class of Bayes rules which will be used in Section 7.

Given $0 \leq k \leq n$, denote $\beta(W^{\theta}_{k-1},W^{\theta}_k)= \exp \{g(W^{\theta}_{k-1},W^{\theta}_k)\}$.
For $0 < p <1$ and $q=1-p$, let
\begin{eqnarray}\label{rnp}
R_{n,p}  &:=& \sum_{k=1}^n \frac{1}{q} \frac{p_n(Y_k,Y_{k+1},\ldots,Y_n;\theta)}
{p_n(Y_k,Y_{k+1},\ldots,Y_n;\theta_0)} \\
&=& \sum_{k=1}^n \frac{1}{q} \beta(W^{\theta}_{n-1},W^{\theta}_n)  \cdots  \beta(W^{\theta}_{k-1},W^{\theta}_k). \nonumber
\end{eqnarray}
%{\rm where}~$Y^{\theta}_{-1}=Y^{\theta}_0.$
%$\pi_{\theta_0}(\cdot), \pi_{\theta_1}(\cdot) \in {\bf M}$.
Note that
%\begin{eqnarray*}
$R_{n+1,p}  =  \beta(W^{\theta}_{n},W^{\theta}_{n+1}) \frac{1}{q}  (1 + R_{n,p}).$
%\end{eqnarray*}
Define
\begin{eqnarray*}
N_{q,b} &=& \inf\{n: R_{n,p} \geq B \}=\inf\{n: R_{n,p} \geq B(W_n) \}, \label{fuh46}\\
H_n(y,w) &=& \mathbb{P}_{\infty} \{ R_{n,p} \leq y|  N_{q,b} > n, W_{n} =w\}, \nonumber \\
\rho(t,y,w) &=& \mathbb{P}_{\infty} \{ R_{n+1,p} \leq y| R_{n,p}=t,  N_{q,b} > n+1, \\
&&~~~~~W_{n+1}=w\}, \nonumber\\
\zeta(t,w,w') &=& \mathbb{P}_{\infty} \{ N_{q,b} > n+1, W_{n+1} \in dw'| R_{n,p}=t,   \\
&&~~~~~N_{q,b} > n,W_n =w\}. \nonumber
\end{eqnarray*}

For a given set of non-negative boundary points $B = \{B(w): w \in {\cal W}\}$
(infinity is not excluded), consider the set $S_B=\{(r,w): w \in {\cal W}, 0 < r < B(w)\}.$
Let ${\cal F}_B$ be the set of distribution functions with support in $S_B$.
For given $H(\cdot,\cdot) \in {\cal F}_B$, let $T_B$ be the transformation on ${\cal F}_B$ defined by
\begin{eqnarray}\label{4.2bay}
&~& T_B H(r,w) 
= \frac{1}{{\mathbb Q}(H)}\int_{w' \in {\cal W}} \int_0^{B(w')} \rho(t,r,w) \zeta(t,w',w) \nonumber\\
&~&~~~~~~~~~~~~~~~~~dH(t,w')\mathbb{P}(w',dw),  
\end{eqnarray}
where
\begin{eqnarray}\label{4.3bay}
{\mathbb Q}(H)=\int_{w,w'\in {\cal W}} \int_0^{B(w')}  \zeta(t,w',w)  dH(t,w')\mathbb{P}(w',dw).
\end{eqnarray}

The following proposition characterizes the behavior of $T_B$.

\begin{proposition}\label{A1}
For each given $B$, we have $T_B H_n=H_{n+1}$. Therefore there associates a set of invariant measures
$\Phi_B$ such that $T_B \phi= \phi$ for all $\phi \in \Phi_B$.
\end{proposition}

The proof of Proposition \ref{A1} is given in the Appendix.\\

By Proposition \ref{A1}, we have that for each $B$ there is an associated set of invariant measures $\Phi_B$, i.e.,
$T_B \phi= \phi$ for all $\phi \in \Phi_B$. Define $\tilde{\phi}$ as
\begin{eqnarray*}
d \tilde{\phi}(y,w) = \frac{\int_{w' \in {\cal W}}(1+py)d\phi(y,w)\mathbb{P}(w,dw')}
{\int_{w,w' \in {\cal W}} \int_{0}^{B(w')} (1+pt)d\phi(t,w)\mathbb{P}(w,dw')}.
\end{eqnarray*}
It is easy to see that if the distribution of $R_{0,p}$ is $\tilde{\phi}$, then the distribution of
$R_{0,p}$ conditional on $\{\omega >0\}$ is $\phi$. Note that $\phi$ depends on $p$.
Let $0 < c < \infty$ and $0 < p <1$ be such that $N_{q,b}$ is the Bayes rule for $B(0,p,c)$.
By using the same argument as that in Theorem 4 of Fuh \cite{Fuh2004b},
we can choose a subsequence $\{T_B^i, p_i, c_i, \phi_i\}$ such that as $i \to \infty, p_i \to 0, c_i \to c^*$ and
$\phi_i$ converges in distribution to a limit $\psi$.

Given the value of the initial state $W_0= \tilde{w}$, the initial $R_0^*(\theta)$ is simulated from
the distribution $\psi$, conditioned on the event $\{W_0=\tilde{w}\}$. Define recursively
\begin{eqnarray}\label{r}
 R_{n+1}^*(\theta)  =  \beta(W^{\theta}_{n},W^{\theta}_{n+1})  (1 + R^*_{n}(\theta)).
\end{eqnarray}
Let
\begin{eqnarray}\label{rf}
  R_n^*(F) = \int_{\theta \in J} R_n^*(\theta) dF(\theta).
\end{eqnarray}
Denote $b=\log B$, and define the weighted Shiryayev-Roberts-Pollak (SRP) rule as
\begin{eqnarray}\label{srp}
 N_{b}^{\psi} := \inf\{n: R_{n}^*(F) \geq B \}=\inf\{n: \log R_{n}^*(F) \geq b \}.
\end{eqnarray}

Note that each one of these detection policies (\ref{fuh46}) and (\ref{srp}) is an ``equalizer rule''  in the sense that
\begin{eqnarray}\label{el}
\mathbb{E}_k ( N_{b}^{\psi} - k+1| N_{b}^\psi \geq k-1)=\mathbb{E}_1 N_{b}^\psi,
\end{eqnarray}
for all $k>1$. The same is true for the case where $\psi$ has atoms on the boundary,
since the randomization law is time independent. Note that the threshold of the Bayes rule (\ref{fuh46}) depends on the current state
of the Markov chain, while the threshold of the SRP rule (\ref{srp}) is a constant. 
By using an argument similar to Lemma 7 of Fuh \cite{Fuh2004b}, we have that the difference 
between these two rules is $o(1)$ as $p \to 0$ and $b \to \infty$.

Next, we will study asymptotic approximations for the average run length in the
weighted SRP detection rule when $w$ is finite.
Since $N_b^\psi$ is an equalizer rule, we only consider the approximation of
${\mathbb E}_1 N_b^\psi$. Given $\theta=\theta^0$ or $\theta \in J$, let $\pi_{\theta}$ denote
the stationary distribution of $\{X_n,n \geq 0\}$ under $P^{\theta}$.
For given ${P}^{\theta_0}$ and ${P}^{\theta}$ and denote
${\mathbb P}^{\theta_0}$ and ${\mathbb P}^{\theta}$ as the induced probabilities,
define the Kullback-Leibler information number
\begin{eqnarray}\label{KL5.9}
 K({P}^{\theta}, {P}^{\theta_0})=K({\mathbb P}^{\theta}, {\mathbb P}^{\theta_0}) = {\mathbb E}^{\theta}
\bigg( \log \frac{\|{\bf P}_{\theta}(Y_1)  \pi_{\theta} \| }
 {\|{\bf P}_{\theta_0}(Y_1)   \pi_{\theta_0} \|}  \bigg).
\end{eqnarray}
By assumption C3, we have
%\begin{equation}
$0<K({P}^{\theta}, {P}^{\theta_0})<\infty.$
%\end{equation}

To derive a second-order approximation for the average run lengths of the weighted SRP rule, we will
apply relevant results from nonlinear Markov renewal theory developed in Section 3 of Fuh \cite{Fuh2004b}.
To this end, we rewrite the stopping time $N_b:=N_b^{\psi}$
(we delete $\psi$ for simplicity) in
the form of a Markov random walk crossing a constant threshold plus
a nonlinear term that is slowly changing.
Note that the stopping time $N_b$ can be written in the following form
\begin{equation}\label{nrn}
N_b =\inf \{n\geq 1: {\mathbb S}_n  + \eta_n \geq b \},~~~b=\log B,
\end{equation}
where for $n \geq 1$,
\begin{eqnarray}\label{sn}
&~& {\mathbb S}_n = {\mathbb S}_n(\theta)  \\
&=& (\theta - \theta_0) \sum_{k=1}^n g(W_{k-1},W_k) - n
(\Lambda(\theta) - \Lambda(\theta_0)), \nonumber
\end{eqnarray}
is a Markov random walk with mean 
${\mathbb E}^\theta {\mathbb S}_1  =
K({\mathbb P}^{\theta},{\mathbb P}^{\theta_0})$, and
\begin{eqnarray}\label{eta}
\eta_n &=& \eta(\theta) \\
& =& \log \int_{\Theta} \frac{r(W_n;\alpha)}{r(W_0;\alpha)}
\exp\bigg\{(\alpha - \theta) \sum_{k=1}^n g(W_{k-1},W_k)  \nonumber \\
&& - n (\Lambda(\alpha) - \Lambda(\theta))\bigg\}  \left\{1+\sum_{k=1}^{n} \exp(-{\mathbb S}_k(\alpha))
 \right\}dF(\alpha). \nonumber
\end{eqnarray}

Suppose there exists $\theta$ such that $\Lambda'(\theta)={\mathbb E}_{\Pi} g(W_0,W_1)$, and
denote $\hat{\theta}_n = \theta(\sum_{k=1}^n g(W_{k-1},W_k)/n)$. Then
$\eta_n$ can be further decomposed as $l_n+ V_n$, where
\begin{eqnarray}
 l_n &=& -\frac{1}{2} \log n, \label{ln} \\
V_n &=& (\hat{\theta}_n - \theta) \sum_{k=1}^n g(W_{k-1},W_k) - n
 (\Lambda(\hat{\theta}_n) - \Lambda(\theta)) \label{vna}  \\
&+& \log n^{1/2} \int_{\Theta} \frac{r(W_n;\alpha)}{r(W_0;\alpha)}
  \exp\bigg\{(\alpha - \hat{\theta}_n) \sum_{k=1}^n g(W_{k-1},W_k)   \nonumber \\
&-& n (\Lambda(\alpha) - \Lambda(\hat{\theta}_n)) \bigg\}  \left\{1+\sum_{k=1}^{n} \exp(-{\mathbb S}_k(\alpha))
 \right\}dF(\alpha) \nonumber \\
&:=& n K({\mathbb P}^{\theta},{\mathbb P}^{\hat{\theta}_n}) + \log u_n (\sum_{k=1}^n g(W_{k-1},W_k)/n). \nonumber
\end{eqnarray}

For $b>0$, define
\begin{equation}\label{nb}
N^*_b=\inf \{n \geq 1: {\mathbb S}_n \geq b\},
\end{equation}
and let $R_b= {\mathbb S}_{N^*_b}-b$ (on $\{N^*_b < \infty\}$) denote the overshoot of the statistic
${\mathbb S}_n$ crossing the threshold $b$ at time $n=N^*_b$. When $b=0$, we denote $N^*_b$ in
($\ref{nb}$) as $N^*_+$. For given $\tilde{w} := (x_0, \pi)\in {\cal W}$, with $x_0$
is the initial state of $X_0$ taken from $\pi(x_0)$, let
\begin{equation}
G(u)=\lim_{b\rightarrow \infty} {\mathbb P}^\theta \{ R_b \leq u |W_0= \tilde{y} \}
\end{equation}
be the limiting distribution of the overshoot. It is known [cf. Theorem 1 of Fuh \cite{Fuh2004a}] that
\begin{eqnarray}
 \lim_{b\rightarrow \infty} {\mathbb E}^{\theta} (R_b|W_0= \tilde{w})
 = \int_0^\infty u dG(u) = \frac{{\mathbb E}_{m_+}^\theta S_{N^*_+}^2}{2 {\mathbb E}_{m_+}^\theta S_{N^*_+}},
\end{eqnarray}
where $m_+:=m_+^\theta$ is defined in the same way as $\pi^\theta_+$ defined in Section 3 of Fuh \cite{Fuh2004b}.
%Let us also define
%\[ \zeta =\lim_{b\rightarrow \infty} {\mathbb E}_1 ( e^{- R_b}|W_0= \tilde{w} ) = \int_0^\infty e^{-y} dG(y). \]
%and
%\begin{equation}C(\rho,D)=\b E_1\log \left\{1+\sum_{i=1}^{\infty}(1-\rho)^ie^{-Z_i}\right\}.\end{equation}

Note that by ($\ref{nrn}$), we have
\begin{eqnarray}\label{4.17a}
 {\mathbb S}_{N_b} = b - \eta_{N_b} + O_b~~~\mbox{ on } \{N_b < \infty\},
\end{eqnarray}
where $O_b= {\mathbb S}_{N_b} + \eta_{N_b}-b$ is the overshoot of ${\mathbb S}_n + \eta_n$
crossing the boundary $b$ at time $N_b$. Taking the expectations on both sides of (\ref{4.17a}), 
and applying Wald's identity for Markov random walks [cf. Corollary 1 of Fuh and Zhang \cite{Fuh2000}], we obtain
\begin{eqnarray}
&~& K({\mathbb P}^{\theta},{\mathbb P}^{\theta_0}) {\mathbb E}^\theta( N_b|W_0= \tilde{w}) \nonumber \\
&~& + \int_{\cal W} \Delta_\theta(w) m^\theta_+(dw) - \Delta_\theta(\tilde{w}) \nonumber  \\
&=& {\mathbb E}^\theta ( {\mathbb S}_{N_b}|W_0= \tilde{w} ) \\
&=& b- {\mathbb E}^\theta
( \eta_{N_b} |W_0= \tilde{w} ) + {\mathbb E}^\theta( O_b|W_0= \tilde{w}), \nonumber
\end{eqnarray}
where $\Delta_\theta : {\cal W} \rightarrow  {\bf R}^d$ solves the Poisson equation
\begin{eqnarray}\label{pois}
\mathbb{E}_{w}^\theta \Delta_\theta(W_1) - \Delta_\theta(w) = \mathbb{E}_{w}^\theta {\mathbb S}_1 -
\mathbb{E}_{m_+}^\theta \mathbb{S}_1
\end{eqnarray}
for almost all $w \in {\cal W}$ with $\mathbb{E}_{m_+}^\theta \Delta_\theta (W_1) = 0$.

The crucial observations are that the sequence $\{V_n, n\geq 1\}$ is slowly changing, and that $V_n$
converges in ${\mathbb P}^\theta$-distribution, as $n \rightarrow \infty$, to the random variable
\begin{eqnarray}\label{4.19}
&~& \tilde{V} \\
&=& \frac{1}{2} \chi_1^2 + \frac{1}{2}\log\frac{2\pi F'(\theta)}
{\Lambda''(\theta)} +  \log \left\{1+\sum_{k=1}^{\infty} \exp (- {\mathbb S}_k(\theta) ) \right\} \nonumber  \\
&&+ \log \bigg\{\frac{\mathbb{E}_{m_+}^{\theta}r(W_{N_+^*};\theta)}{r(W_0; \theta)}\bigg\}, \nonumber
\end{eqnarray}
where $\chi_1^2$ denotes a random variable having the chi-squared distribution with one degree of
freedom.

Denote $\gamma_\theta = \log \left\{1+ \sum_{k=1}^{\infty} \exp (- {\mathbb S}_k (\theta) ) \right\}$, we will show in
Section 6 that ${\mathbb E}_{m_+}^\theta \gamma_\theta < \infty.$
Here the expectation ${\mathbb E}_{m_+}^\theta$ is taken under
$\omega=1$ and the initial distribution of $Y_0$ is $m_+$, we omit $1$ for simplicity.
An important consequence of the slowly changing property is that, under mild conditions,
the limiting distribution of the overshoot of a Markov random walk over a fixed threshold does not change by the
addition of a slowly changing nonlinear term [cf. Theorem 1 in Section 3 of Fuh \cite{Fuh2004b}].
More importantly, nonlinear Markov renewal theory allows us to obtain an asymptotically accurate approximation
for ${\mathbb E}N_b$, that takes the overshoot into account. Now we can characterize the constant $C(\theta)$
appeared in Theorem \ref{thm2},
\begin{eqnarray}\label{aec}
&C(\theta) 
 = \frac{{\mathbb E}_{m_+}^\theta S_{N^*_+}^2}{2 {\mathbb E}_{m_+}^\theta S_{N^*_+}}
- {\mathbb E}_{m_+}^\theta  \gamma_\theta  - \frac{1}{2}\log\frac{2\pi F'(\theta)}{\Lambda''(\theta)}  - \frac{1}{2} \\
 & (\int_{\cal W} \Delta(w) m^\theta_+(dw) - \Delta(\tilde{w})) 
- \log \bigg\{ \frac{\mathbb{E}_{m_+}^{\theta}r(Y_{N_+^*};\theta)}{r(Y_0; \theta)}\bigg\}. \nonumber
\end{eqnarray}

%\cdf{double check the constant}

When $\theta=\theta_1$ is known, we have the following approximation of the
average run length. Since the proof is similar to that of Theorem 6 in Fuh \cite{Fuh2004b}, we
will not repeat it here.
\begin{proposition}
Let $Y_1,\ldots,Y_n$ be a sequence of random
variables from a state space model $\{Y_n, n \geq 1\}$ satisfying conditions {\rm C1-C4}.
Assume that ${S}_1$ is nonarithmetic with respect to
${P}_{\infty}$ and ${P}_1$.
%If $0< K({\mathbb P}^{\theta_1},{\mathbb P}^{\theta_0}) < \infty$,
%$0 < K({\mathbb P}^{\theta_0},{\mathbb P}^{\theta_1}) < \infty$,
%and ${\mathbb E}_1|{\mathbb S}_1|^2 < \infty$,
Then for $\tilde{w} \in {\cal W}$, as $b \rightarrow \infty$
\begin{eqnarray}
&~& {\mathbb E}_1 (N_b|W_0=\tilde{w}) \\
&=& \frac{1}{K({\mathbb P}^{\theta_1},{\mathbb P}^{\theta_0})}
\bigg(b- {\mathbb E}_{m_+} \gamma + \frac{{\mathbb E}_{m_+} S_{N^*_+}^2}{2 {\mathbb E}_{m_+} S_{N^*_+}} \nonumber \\
&~&~~~- \int_{\cal W} \Delta(w) m_+(dw) +\Delta(\tilde{w})\bigg) + o(1). \nonumber
\end{eqnarray}
\end{proposition}

\def\theequation{6.\arabic{equation}}
\setcounter{equation}{0}
\section{Proof of Theorem \ref{thm2}}\label{sec6}

To prove Theorem \ref{thm2}, without loss of generality, we assume that
$J=\Theta = [\theta_0,\theta_1] \subset {\bf R}$ and $\theta \geq 0$.
Note that the proof of $(\ref{ae})$ rests on the nonlinear Markov renewal theory from Theorem 3 and Corollary 1 in
Fuh \cite{Fuh2004b}. Indeed, by $(\ref{nrn})$, the stopping time $N_b^{\psi}$ is based on the
thresholding of the sum of the Markov random walk ${\mathbb S}_n$ and the nonlinear term $\eta_n$.
From (\ref{ln}) and (\ref{vna}), we
have $\eta_n=l_n + V_n$, with $l_n=-(1/2) \log n$. It is easy to see that
$\lim_{n \to \infty} \max_{0 \leq j \leq \sqrt{n}} | - (1/2) \log (n+j) + (1/2) \log n | = 0.$
In order to apply Theorem 3 and Corollary 1 in Fuh \cite{Fuh2004b},
we need to check the validity of the conditions which are stated in the following lemmas, respectively.
Relation (\ref{ae}) will then follow by specialization. Note that although the nonlinear Markov renewal
theory developed in Fuh \cite{Fuh2004b} is under the condition of $w$-uniformly ergodic, it can be generalized
to the norm in Definition \ref{def3}. A heuristic explanation of this result can be described as follows:
we first investigate the difference
between a stopping time crossing nonlinear boundaries and a stopping time crossing linear boundaries with varying
drift, then derive nonlinear Markov renewal theory directly from parallel results in the linear case with varying drift
via the uniform integrabilities and the weak convergence of the overshoot.
Because the uniform Markov renewal theory developed in Fuh \cite{Fuh2004a} is under a general norm, therefore the extension of the proofs
in Fuh \cite{Fuh2004b} is straightforward. The details are omitted.

In the proof of the following lemmas, we will assume the conditions of Theorem \ref{thm2} hold.
We first consider the case that $F$ is concentrated on $[\theta_0, \theta_1]$,
where $ 0 < \theta_0 <\theta < \theta_1 < \infty$
are such that $\alpha \Lambda'(\alpha)- \Lambda(\alpha) > 0$ for $\theta_0 \leq \alpha \leq \theta_1$
and $F$ has a derivative $F'$ which is positive and continuous on $[\theta_0,\theta_1]$.
The probability $\mathbb{P}_1$ and expectation $\mathbb{E}_1$
in this section are taken under $Y_0=\tilde{y}$, and we omit it for simplicity.

\begin{lemma}\label{61}
Under assumptions of {\rm Theorem \ref{thm2}},
${\mathbb S}_1 = {\mathbb S}_1(\theta) = (\theta - \theta_0) g(Y_{0},Y_1) -
(\Lambda(\theta) - \Lambda(\theta_0))$ has a nonarithmetic distribution.
\end{lemma}

PROOF.~Suppose the ${\mathbb S}_1(\theta)$ has an arithmetic distribution for some
$\theta \neq \theta_0$, say $\theta=\theta^*$, and let $d_1$ be the span of ${\mathbb S}_1(\theta^*)$.
Then $g(Y_{0},Y_1)$ must take values of the form
%\begin{eqnarray*}
$(\frac{1}{\theta^*-\theta_0})\{kd_1+[\Lambda(\theta^*)-\Lambda(\theta_0)]\}=d k+\gamma,~ \mbox{say},$
%\end{eqnarray*}
where $k=0,\pm 1,\pm 2,\ldots.$ Moreover, since $g(Y_{0},Y_1)$ is assumed to have a nondegenerate distribution,
there are $k_1\neq k_2$ for which ${\mathbb P}_y \{d k_1+\gamma\}>0< {\mathbb P}_y \{d k_2+\gamma\}$ for all
$y \in {\cal Y}$. Now suppose that
${\mathbb S}_1(\theta)$ has an arithmetic distribution for some $\theta$ with $\theta_0\neq \theta\neq \theta^*$
and let $d(\theta)>0$ denote the span of ${\mathbb S}_1(\theta)$. Then there are $j_1$ and $j_2$ for
which $j_1 \neq j_2$ and
%\begin{eqnarray*}
$j_i d(\theta)=(\theta-\theta_0)(d k_i+\gamma)-[\Lambda(\theta)-\Lambda(\theta_0)],~i=1,2.$
%\end{eqnarray*}
Therefore
\begin{equation}\label{6.9}
\frac{\Lambda(\theta)-\Lambda(\theta_0)}{\theta-\theta_0}=\gamma+d(\frac{k_1 j_2-k_2j_1}{j_2-j_1}).
\end{equation}
Thus, the set of $\theta$ for which $\theta_0\neq \theta\neq \theta^*$ and ${\mathbb S}_1(\theta)$
has an arithmetic distribution is
contained in the set of $\theta$ for which (\ref{6.9}) holds for some $j_1\neq j_2$;
the latter set is countable, since $\Lambda$ is convex.~~~~~~~$\Box$

\begin{lemma}
Under assumptions of {\rm Theorem \ref{thm2}}, we have
\begin{eqnarray}
 \sum_{n=1}^\infty {\mathbb P}_1 \{V_n \leq -\varepsilon n\} < \infty \mbox{ for some }
0< \varepsilon < K(\mathbb{P}^{\theta_1},\mathbb{P}^{\theta_0}). \label{82}
\end{eqnarray}
\end{lemma}
Condition (\ref{82}) holds trivially because $r(y;\theta)$ is uniformly positive and hence $V_n \geq 0$.

\begin{lemma}
Under assumptions of {\rm Theorem \ref{thm2}}, then
\begin{eqnarray}
 \max_{0\leq l\leq n}|V_{n+l}|,~ n\geq 1,~ \mbox{are}~{\mathbb P}_1-uniformly~integrable. \label{83}
\end{eqnarray}
\end{lemma}

PROOF.~To show (\ref{83}) holds, we first prove
\begin{eqnarray}
\max_{0\leq l \leq n} (\hat{\theta}_{n+l} - \theta) \sum_{k=1}^{n+l} g(Y_{k-1},Y_k) - (n+l)
 (\Lambda(\hat{\theta}_{n+l}) - \Lambda(\theta)) \label{83a}
\end{eqnarray}
are ${\mathbb P}_1$-uniformly integrable, where $\hat{\theta}_{n+l}$ is the maximum likelihood
estimator of $\theta$.
Note that on the event $A_n$ of $ |(1/n )\sum_{k=1}^n g(Y_{k-1},Y_k) - \Lambda'(\theta)| < \varepsilon$ for some
$\varepsilon > 0$, we have for all $n \geq 1$,
\begin{eqnarray*}
&~& (\hat{\theta}_n - \theta) \sum_{k=1}^n g(Y_{k-1},Y_k) - n
 (\Lambda(\hat{\theta}_n) - \Lambda(\theta)) \\
&\leq&  B n \bigg(\frac{1}{n} \sum_{k=1}^n g(Y_{k-1},Y_k)
- \Lambda'(\theta)\bigg)^2
\end{eqnarray*}
on $A_n$, for some constant $B$. Therefore,
\begin{eqnarray}
&~& {\mathbb P}_1 \bigg\{ \max_{0\leq l \leq n} \{(\hat{\theta}_{n+l} - \theta) \sum_{k=1}^{n+l} g(Y_{k-1},Y_k)  \nonumber \\ 
&~&~~~~~-(n+l)
 (\Lambda(\hat{\theta}_{n+l}) - \Lambda(\theta))\} > a \bigg\}  \\
&\leq&  {\mathbb P}_1 \bigg\{ \max_{0\leq l \leq 2n} l | \frac{1}{l} \sum_{k=1}^l g(Y_{k-1},Y_k)
- \Lambda'(\theta)| > \sqrt{Bna} \bigg\}. \nonumber \label{quick11}
\end{eqnarray}
Since conditions of Theorem 2 imply that conditions of Theorem 2 in Fuh and Zhang \cite{Fuh2000} hold, we have that
for all $\varepsilon>0$ and $r\geq 0$
\begin{eqnarray}\label{quick1}
\sum_{n=1}^\infty n^{r-1} {\mathbb P}_1 \left\{\max_{1\leq l\leq n}({\mathbb S}_l -
 \Lambda(\theta) l ) \geq \varepsilon n\right\} < \infty.
\end{eqnarray}
Hence ${\rm (\ref{quick11})} \leq C a^{-r}$, for some $C > 0$ and $r > 1$. This imply (\ref{83a}) hold.

Denote
\begin{eqnarray}
&~& H_n \\ 
&=& n^{1/2} \int_{\Theta} \frac{r(Y_n;\alpha)}{r(Y_0;\alpha)}
   \exp\bigg\{(\alpha - \hat{\theta}_n) \sum_{k=1}^n g(Y_{k-1},Y_k) \nonumber \\
&~& - n (\Lambda(\alpha) - \Lambda(\hat{\theta}_n)) \bigg\} \nonumber
\times \left\{1+\sum_{k=1}^{n-1} \exp(-{\mathbb S}_k(\alpha)) \right\}dF(\alpha). \label{6.37}
\end{eqnarray}
To complete the proof, we need to show that $\max_{0\leq l \leq n} H_{n+l}$ are
${\mathbb P}_1$-uniformly integrable. First, we note that
$(\hat{\theta}_{n+l} - \theta) \sum_{k=1}^{n+l} g(Y_{k-1},Y_k) - (n+l)
 (\Lambda(\hat{\theta}_{n+l}) - \Lambda(\theta))$ are uniformly bounded on $A_n$ and
$0 < r(y;\theta) < \infty$ uniformly for $y \in {\cal Y}$ by Theorem 4 of Chan and Lai \cite{Chan2003}.

To analyze the term appeared in (\ref{6.37}), denote
$W^{n}_{\alpha}= 1 + \sum_{k=1}^{n-1}\exp (- {\mathbb S}_{k}(\alpha) ),$ for $\theta_0 \leq \alpha \leq \theta_1$.
Note that $W^{n}_{\alpha}$ converges ${\mathbb P}_1^{\theta}$-a.s. as $n\rightarrow \infty$ to
a random variable $W_{\alpha}^{\theta} :=1 + \sum_{k=1}^{\infty} \exp (- {\mathbb S}_{k}(\alpha))$. Since
\begin{eqnarray*} 
&~& \sum_{n=m}^{\infty}(W^{n+1}_{\alpha}-W^{n}_{\alpha}) \\
& =& \sum_{n=m}^{\infty} \exp \left \{- \left[
\alpha  \sum_{k=1}^n g(X_{k-1},X_k) - n \Lambda(\alpha)\right ] \right \}
\rightarrow_{m\rightarrow \infty}0~
\end{eqnarray*}
${\mathbb P}_1^{\theta}-a.s.$, uniformly in $\alpha \in [\theta_0, \theta_1]$, it follows that
$W_{\alpha}^{\theta}$ is ${\mathbb P}_1^{\theta}-a.s.$ continuous in $\alpha \in [\theta_0, \theta_1]$.

Next we will show, which is more than it suffices, that there exists a constant $a >0$ such that
\begin{eqnarray}\label{p83}
{\mathbb E}_1^{\theta}\left (\int_{\theta_0}^{\theta_1} \left\{1+\sum_{k=1}^{\infty} \exp(-{\mathbb S}_k(\alpha)) \right\}dF(\alpha) \right)^a < \infty.
\end{eqnarray}

For given $\varepsilon >0$, let $\Gamma=\min\{n||\sum_{k=1}^m g(Y_{k-1},Y_k)/m - \Lambda'(\theta) |
\leq \varepsilon~{\rm for~all}~m\geq n\}$.
Suppose that $\varepsilon$ is chosen small enough so that there exists $\beta>0$ such that
$ {\mathbb S}_n(\alpha) \geq \beta n $ if $n \geq \Gamma$ for all $\theta_0 \leq \alpha \leq \theta_1$.
There exists a constant $\eta >0$ such that
$ | \Lambda(\theta-\alpha)+ \Lambda(\alpha)- \Lambda(\theta)| < \eta$ for all
$\theta_0 \leq \alpha \leq \theta_1.$ By using the large deviation result for Markov random walks (cf. Ney and Nummelin \cite{Ney1987}) we can choose a constant $\delta>0$ such that ${\mathbb P}_1^{\theta}(\Gamma = \lambda)
\leq \exp\{- \delta \lambda \}$. Furthermore we choose $1>a>0$ such that $a \eta - \delta(1-a)<0$. Now
\begin{eqnarray}
&~& \int_{\theta_0}^{\theta_1} W_{\alpha}^{\theta} dF(\alpha) \label{6a} \\
&=& \int_{\theta_0}^{\theta_1} \left( 1+ \sum_{k=1}^{\Gamma-1} \exp (- {\mathbb S}_{k}(\alpha) )+
  \sum_{k= \Gamma}^{\infty} \exp (-  {\mathbb S}_{k}(\alpha)) \right) dF(\alpha) \nonumber \\
 &\leq& \int_{\theta_0}^{\theta_1} \left( 1+ \sum_{k=1}^{\Gamma-1} \exp (- {\mathbb S}_{k}(\alpha) )+
  \frac{1}{1- e^{- \beta}} \right) dF(\alpha). \nonumber 
\end{eqnarray}
To evaluate the second term in the integrand of (\ref{6a}), we have
\begin{eqnarray*}
&~& {\mathbb E}_1^{\theta} \left( \int_{\theta_0}^{\theta_1} \sum_{k=1}^{b-1} \exp (- {\mathbb S}_{k}(\alpha))  dF(\alpha)\bigg|\Gamma = b \right) \\
& \leq& \displaystyle \frac{{\mathbb E}_1^{\theta}
\int_{\theta_0}^{\theta_1} \sum_{k=1}^{b-1} \exp (- {\mathbb S}_{k}(\alpha)) dF(\alpha)}
{{\mathbb P}_1^{\theta} (\Gamma=b)} \\
&=&\displaystyle \frac{1}{{\mathbb P}_1^{\theta}(\Gamma=b)}
  \int_{\theta_0}^{\theta_1}\sum_{k=1}^{b-1} \{ e^{[\Lambda(\theta-\alpha)+\Lambda(\alpha)-\Lambda(\theta)]k} + O(\rho^k)\}dF(\alpha) \\
&\leq& \displaystyle \frac{1}{{\mathbb P}_1^{\theta}(\Gamma=b)}\bigg( \frac{1}{\eta}e^{\eta b} + \frac{\rho^2}{1-\rho}
\bigg),
\end{eqnarray*}
where $0< \rho < 1$. By Jensen's inequality,
\begin{eqnarray}\label{egamma}
&~& {\mathbb E}_1^{\theta}\left (\int_{\theta_0}^{\theta_1} W_{\alpha}^{\theta}dF(\alpha) \right)^a \\
&=& {\mathbb E}_1^{\theta}\bigg({\mathbb E}_1^{\theta}\bigg[ \left (\int_{\theta_0}^{\theta_1} W_{\alpha}^{\theta} dF(\alpha) \right)^a \bigg| \Gamma \bigg] \bigg) \nonumber \\
& \leq& \sum_{b=1}^{\infty} \left( \displaystyle \frac{1}{{\mathbb P}_1^{\theta} (\Gamma=b)}\bigg(
 \frac{1}{\eta}e^{\eta b}+ \frac{\rho^2}{1-\rho} \bigg) + \frac{2-e^{-\beta}}{1-e^{-\beta}} \right)^a \nonumber \\
&~&~~~{\mathbb P}_1^{\theta}(\Gamma=b). \nonumber
\end{eqnarray}
The inequality (\ref{p83}) now follows because there exist constants $C_1,\cdots,C_a$ such that
\begin{eqnarray*}
&~& \displaystyle \sum_{b=1}^{\infty} \left( \frac{1}{{\mathbb P}_1^{\theta}(\Gamma=b)}
\bigg(\frac{e^{\eta b}}{\eta}+ \frac{\rho^2}{1-\rho}\bigg) \right)^a {\mathbb P}_1^{\theta}(\Gamma=b) \\
&=& \displaystyle \sum_{b=1}^{\infty} \bigg(\frac{e^{ \eta b}}{\eta}+\frac{\rho^2}{1-\rho}\bigg)^a
[ {\mathbb P}_1^{\theta}(\Gamma=b)]^{1-a} \\
&\leq& \displaystyle \frac{1}{\eta^a} \sum_{b=1}^{\infty} e^{b (a \eta - \delta(1-a))} + C_1 \displaystyle
\frac{1}{\eta^a} \sum_{b=1}^{\infty} e^{b ((a-1) \eta - \delta(1-a))} + \cdots \\
&~&~~~+ C_a
\displaystyle \frac{1}{\eta^a} \sum_{b=1}^{\infty} e^{b ( \eta - \delta(1-a))}
<\infty. ~~~\Box
\end{eqnarray*}

\begin{lemma}
Let $V_n$ be defined in {\rm (\ref{vna})} and $\tilde{V}$ be defined in {\rm (\ref{4.19})}. Then under assumptions of {\rm Theorem \ref{thm2}}, we have
\begin{eqnarray}
&~& V_n \longrightarrow_{\scriptstyle n\rightarrow \infty} \tilde{V}~~~in~{\mathbb P}_{1}\mbox{-distribution} \\
&~&~\mbox{and}~~~{\mathbb E}_{1} V_n \longrightarrow_{\scriptstyle n\rightarrow \infty} {\mathbb E}_{1} \tilde{V}.
\nonumber \label{81}
\end{eqnarray}
\end{lemma}

PROOF. Let $A_n = \{ |\sum_{k=1}^n g(Y_{k-1},Y_k)/n - \Lambda'(\theta)| < \varepsilon\}$. Then by using a result of
large deviations in Markov random walks [cf. Ney and Nummelin \cite{Ney1987}],
there exists a $\delta > 0$ such that ${\mathbb P}_1 \{ A_n^{c}\} \leq \delta$.
Let $\theta$ be defined such that $\Lambda'(\theta)={\mathbb E}_{\Pi} g(Y_0,Y_1)$.
Under the event $A_n$, the maximum likelihood estimate $\hat{\theta}_n=\theta(\sum_{k=1}^n g(Y_{k-1},Y_k)/n)$
is well defined.
Recall $\eta_n = \l_n + V_n,$ where $\eta_n$ is defined in (\ref{eta}), $\l_n = (-1/2) \log n$, and
\begin{eqnarray}
&~& V_n \label{vb}\\
&=& (\hat{\theta}_n - \theta) \sum_{k=1}^n g(Y_{k-1},Y_k) - n
 (\Lambda(\hat{\theta}_n) - \Lambda(\theta)) \nonumber \\
&~& + \log n^{1/2} \int_{\Theta} \frac{r(Y_n;\alpha)}{r(Y_0;\alpha)}
   \nonumber\\
&~& \times  \exp\bigg\{(\alpha - \hat{\theta}_n) \sum_{k=1}^n g(Y_{k-1},Y_k) - n
(\Lambda(\alpha) - \Lambda(\hat{\theta}_n)) \bigg\} \nonumber\\ &~&~~ \left\{1+\sum_{k=1}^{n} \exp(-{\mathbb S}_k(\alpha))
 \right\}dF(\alpha) \nonumber \\
&:=& n K({\mathbb P}^{\theta},{\mathbb P}^{\hat{\theta}_n}) + \log u_n (\sum_{k=1}^n
g(Y_{k-1},Y_k)/n). \nonumber
\end{eqnarray}

We first analyze the second term in (\ref{vb}) and show that for any $\theta \in \Theta$
\begin{eqnarray}\label{6.69}
&~& log u_n (\sum_{k=1}^n g(Y_{k-1},Y_k)/n) \\
&\longrightarrow& \frac{1}{2}\log\frac{2\pi F'(\theta)}
{\Lambda''(\theta)} +  \log \left\{1+\sum_{k=1}^{\infty} \exp (- {\mathbb S}_k(\theta) ) \right\} \nonumber \\
&~&~+ \log \bigg\{\frac{\mathbb{E}_{m_+}^{\theta}r(Y_{N_+^*};\theta)}{r(Y_0; \theta)}\bigg\}, \nonumber
\end{eqnarray}
${\mathbb P}_1^\theta$-a.s. as $n \to \infty$.

To complete the proof of (\ref{6.69}).
By (\ref{p83}), and $0 < r(y;\alpha) < \infty$ uniformly for $y \in {\cal Y}$ for all $\alpha \in \Theta$
via Theorem 4 of Chan and Lai \cite{Chan2003}, we need only to show that
\begin{eqnarray}\label{6.70}
&~& \log n^{1/2} \int_{\Theta} \exp\bigg\{(\alpha - \hat{\theta}_n) \sum_{k=1}^n g(Y_{k-1},Y_k) - n
 (\Lambda(\alpha) \nonumber \\
&~&~  - \Lambda(\hat{\theta}_n)) \bigg\} dF(\alpha) \\
&\longrightarrow& \frac{1}{2}\log\frac{2\pi F'(\theta)} {\Lambda''(\theta)}~~~{\mathbb P}_1^\theta-a.s.~as~
n \to \infty, \nonumber
\end{eqnarray}

For given $\alpha \in \Theta$ and $y \in {\bf R}$, let 
$$H(\alpha,y)=(\Lambda(\alpha) - \Lambda(\hat{\theta}_n)) - (\alpha - \hat{\theta}_n)\frac{1}{n} \\ \sum_{k=1}^n g(Y_{k-1},Y_k).$$  Then
\begin{eqnarray}\label{6.71}
{\rm (\ref{6.70})} = \int_\Theta \exp[-n H(\alpha,y)] dF(\alpha),~~~~~y\in {\bf R}.
\end{eqnarray}
Observe that $H(\alpha,y)$ is convex in $\Theta$ for fixed $y \in {\bf R}$, since $\Lambda$ is convex.
Moreover, for fixed $y$, $H(\alpha,y)=\frac{1}{2}\Lambda''(\alpha^\ast)(\alpha-\hat{\theta}_n)^2$,
where $\alpha^\ast=\alpha^\ast(\alpha,y)$ is an intermediate point between $\alpha$ and $\hat{\theta}_n$.
Let $K$ be any compact subinterval of ${\bf R}$. Then there are
a $\sigma > 0$ and a compact $J \subset \Theta$ for which
$[\hat{\theta}_n - \delta,\hat{\theta}_n+\delta] \subset J$ for all
$y\in K$ and, since $\Lambda''$ is positive and continuous, there is an $\varepsilon>0$ for which
$\Lambda''(\alpha^\ast) \geq \varepsilon$ for $|\alpha-\hat{\theta}_n| \leq \delta$ and $y \in K$.
In particular, $H(\alpha,y) \geq \frac{1}{2}\varepsilon(\alpha-\hat{\theta}_n)^2$ for
$|\alpha-\hat{\theta}_n| \leq \delta$ and $y \in K$.
Since $H$ is convex in $\alpha$ for fixed $y$, it follows that $H(\alpha,y)\geq \frac{1}{2}\varepsilon \delta^2$
for $|\alpha-\hat{\theta}_n| \geq \delta$ and $y \in K$ and, consequently, that
\begin{eqnarray}\label{6.72}
\int_{|\alpha-\hat{\theta}_n| \geq \delta} e^{-nH} dF(\alpha) \leq e^{-\varepsilon\delta^2n/2},~y\in K,~n\geq1.
\end{eqnarray}
%Next, let $y_n\in K$,$n \geq 1$, be a convergent sequence, say, $y_n\rightarrow y\in K$ as $n \rightarrow \infty$,
%let $\hat{\omega}_n = \hat{\omega}(y_n)$, $n\geq 1$ and $\hat{\omega} = \hat{\omega}(y)$, and consider the integral over
%$|\omega-\hat{\omega}|<\delta$.
Next, consider the change of variables $\hat{\theta}_n= \hat{\theta}_n+ n^{-1/2}\alpha$ shows that
\begin{eqnarray}\label{6.73}
&~& \sqrt{n}\int_{|\alpha-\hat{\theta}_n|<\delta}\exp\{-nH\}dF(\alpha) \\
&=& \int_{-\delta\sqrt{n}}^{\delta\sqrt{n}}
\exp[-\frac{1}{2} \Lambda''(\alpha_n^\ast) \alpha^2] F'(\hat{\theta}_n+\frac{\alpha}{\sqrt{n}})d\alpha, \nonumber
\end{eqnarray}
where $\alpha_n^\ast=\alpha^\ast(\hat{\theta}_n+n^{-1/2}\alpha,y_n)$, $n\geq 1$.
As $n \rightarrow \infty$, the integrand on the right side of (\ref{6.73}) converges to
$\exp[-\frac{1}{2}\Lambda''(\theta) \alpha^2]F'(\theta)$; and the integrand
is dominated by $C \exp(-\frac{1}{2}\varepsilon \alpha^2)$ for some $C$. So, the right hand side of (\ref{6.73})
converges to
\begin{eqnarray}\label{6.74}
\int_{-\infty}^{\infty} \exp[-\frac{1}{2}\Lambda''(\theta)\alpha^2]F'(\theta)d\alpha
= \sqrt{\frac{2\pi}{\Lambda''(\theta)}} F'(\theta)
\end{eqnarray}
by the dominated convergence theorem.

Finally, using Theorem 17.2.2 of Meyn and Tweedie \cite{Meyn2009}, we have  as
$n \to \infty,$
\begin{eqnarray}\label{6.75}
(\hat{\theta}_n - \theta) \sum_{k=1}^n g(Y_{k-1},Y_k) - n
 (\Lambda(\hat{\theta}_n) - \Lambda(\theta)) \longrightarrow \chi_1^2,
\end{eqnarray}
where $\chi^2_1$ is a random variable with chi-squared distribution with one degree of freedom.

Combining (\ref{6.69}) and (\ref{6.75}), we have the proof.~~~~~~~$\Box$

\begin{lemma}\label{6l}
Under assumptions of {\rm Theorem \ref{thm2}}, we have for some $0<\varepsilon<1$,
\begin{eqnarray}
&~& \lim_{b \rightarrow \infty} b~{\mathbb P}_1 \bigg\{ N_b^{\psi} \leq \frac{\varepsilon b}{K(\mathbb{P}^{\theta},
\mathbb{P}^{\theta_0})}\bigg\}=0. \label{84}
\end{eqnarray}
\end{lemma}

PROOF.~By using ${\mathbb E}_1 g(Y_0,Y_1) > 0$, and $0 < K(\mathbb{P}^{\theta},\mathbb{P}^{\theta_0})< \infty$, we will prove that
\begin{eqnarray}\label{p84}
{\mathbb P}_1 \bigg\{N_b^{\psi} < \frac{(1-\varepsilon)b}{K(\mathbb{P}^{\theta},
\mathbb{P}^{\theta_0})} \bigg\}
\leq e^{- y_\varepsilon b} + \alpha_1(\varepsilon,b),
\end{eqnarray}
where $y_\varepsilon>0$ for all $\varepsilon>0$, and
\begin{eqnarray}\label{alpha1}
\alpha_1(\varepsilon,b)&=&{\mathbb P}_1 \left\{ \max_{1 \leq n<K_{\varepsilon,b}} {\mathbb S}_n
\geq (1+\varepsilon)(1- \varepsilon) b \right\},\\
 K_{\varepsilon,b}&=&\frac{(1-\varepsilon)b}{K(\mathbb{P}^{\theta},\mathbb{P}^{\theta_0})}. \nonumber
\end{eqnarray}
If $(\ref{p84})$ is correct, then the first term on the right hand side of $(\ref{p84})$
is $o(1/b)$ as $b \rightarrow \infty$.
All it remains to do is to show that $\alpha_1(\varepsilon,b)$ in $(\ref{alpha1})$ is $o(1/b)$.

Note that Theorem \ref{thm2} implies conditions of Theorem 2 in Fuh and Zhang \cite{Fuh2000} hold. Hence
for all $\varepsilon>0$ and $r\geq 0$
\begin{eqnarray}\label{quick}
\sum_{n=1}^\infty n^{r-1} {\mathbb P}_1 \left\{\max_{1\leq k\leq n}({\mathbb S}_k-
K(\mathbb{P}^{\theta},\mathbb{P}^{\theta_0}) k ) \geq \varepsilon n\right\} < \infty,
\end{eqnarray}
whenever ${\mathbb E}_1 |{\mathbb S}_1|^2 < \infty$ and ${\mathbb E}_1 [({\mathbb S}_1-
K(\mathbb{P}^{\theta},\mathbb{P}^{\theta_0}))^+]^{r+1}<\infty$. Recall that
under conditions of Theorem \ref{thm2}, ${\mathbb E}_1 |{\mathbb S}_1|^2<\infty$,
and hence, the sum on the left hand side of the inequality $(\ref{quick})$ is finite for $r=1$ and all $\varepsilon>0$,
which implies that the summand should be $o(1/n)$. Since
$\alpha_1(\varepsilon,b)\leq {\mathbb P}_1 \left\{\max_{ n< K_{\varepsilon,b}}({\mathbb S}_n-
K(\mathbb{P}^{\theta},\mathbb{P}^{\theta_0}) n)\geq \varepsilon(1- \varepsilon)b \right\},$
it follows that $\alpha_1(\varepsilon,b)=o(1/b)$.

Next, we need to prove $(\ref{p84})$. We only consider the case that $\theta_0 <\theta$,
as the other case can be done by using a similar way.
Denote $\mathbb{S}^k_n = \log LR_n^k$, and let $N=N_b^{\psi}$ for simplicity.
Let $I\{\cdot\}$ be the indicator function. Recall from (\ref{et}), we have
\begin{eqnarray*}
&~& {\mathbb P}^{\theta_0}(y,dz) \\
&=& \frac{r(z;\theta_0)}{r(y;\theta_0)}\frac{r(y;\theta)}{r(z;\theta)}
     e^{-(\Lambda(\theta_0)- \Lambda(\theta)) + (\theta_0 - \theta) g(Y_0, Y_1)} {\mathbb P}^{\theta}(y,dz).
\end{eqnarray*}
By Proposition 1, for all $\theta \in \Theta$ $0< r(z;\theta) < \infty$ uniformly for $z \in {\cal Y}$.
For any $C >0$, by using a change of measure argument, we have
\begin{eqnarray*}
&~& {\mathbb P}_{\infty} \bigg\{N < (1-\varepsilon)b K(\mathbb{P}^{\theta}, \mathbb{P}^{\theta_0})^{-1} \bigg\} \\
&=& {\mathbb E}_1 \bigg\{I\{ N < K_{\varepsilon,b}\} \frac{r(Y_N;\theta_0)}{r(Y_k;\theta_0)}\frac{r(Y_k;\theta)}
{r(Y_N;\theta)} e^{-(\Lambda(\theta_0)- \Lambda(\theta)) + (\theta_0 - \theta) \mathbb{S}^k_{N}} \bigg\}  \\
&\geq& K {\mathbb E}_1 \bigg\{I \{N < K_{\varepsilon,b},~ \mathbb{S}^k_N < C\} \exp(- k \mathbb{S}^k_N) \bigg\} \\
&\geq& e^{-k C} {\mathbb P}_1 \bigg\{ N < K_{\varepsilon,b},~\max_{n < K_{\varepsilon,b}} \mathbb{S}^k_{n} < C \bigg\}  \\
&\geq& e^{-kC} \bigg[ {\mathbb P}_1 \bigg\{ N < K_{\varepsilon,b} \bigg\} - {\mathbb P}_1
\bigg\{ \max_{n < K_{\varepsilon,b}}\mathbb{S}^k_{n} \geq C \bigg\} \bigg],
\end{eqnarray*}
where $K>0$ is a constant such that $|\frac{r(Y_N;\theta_0)}{r(Y_k;\theta_0)}\frac{r(Y_k;\theta)}
{r(Y_N;\theta)}| > K$, and $k=\theta-\theta_0 > 0$. Choosing $kC \leq (1+\varepsilon)(1-\varepsilon)b$, then, we have
\begin{eqnarray}\label{73pp}
&~& {\mathbb P}_1 \bigg\{N < \frac{(1-\varepsilon)b}{K(\mathbb{P}^{\theta},\mathbb{P}^{\theta_0})} \bigg\} \\
&\leq& e^{kC} {\mathbb P}_{\infty} \bigg\{N < (1-\varepsilon)b K(\mathbb{P}^{\theta},  \mathbb{P}^{\theta_0})^{-1} \bigg\}  + \alpha_1(\varepsilon,b). \nonumber
\end{eqnarray}

Recall that $R_n^*(F)$ is defined in $(\ref{rf})$. Note that under the condition of
$0 < K(\mathbb{P}^{\theta},\mathbb{P}^{\theta_0})< \infty$, we have
%\begin{eqnarray*}
$ {\mathbb P}_{\infty} \big\{N <  K_{\varepsilon,b} \big\}
= \sum_{i=1}^{ K_{\varepsilon,b}} {\mathbb P}_{\infty} \big\{ R^*_i(F) > B \big\}
\leq \sum_{i=1}^{ K_{\varepsilon,b}} \frac{i}{B}
\leq \frac{(\log B)^2}{( K(\mathbb{P}^{\theta},\mathbb{P}^{\theta_0}))^2 B}.$
%\frac{1}{B}{\mathbb E}_{\infty} e^{- (\mathbb{S}_N - b)} (1 + o(1))~~~{\rm as~}B \to \infty.
%\end{eqnarray*}
By letting a suitable $kC$, we have the first term on the right hand side
of $(\ref{73pp})$ $\leq e^{- y_\varepsilon b}$, for some $y_{\varepsilon} > 0$, and get the proof of $(\ref{p84})$.
~~~$\Box$ \\

Next we consider the case that $F$ is a measure on the real line. Assume there exist constants
$0< c< K(\mathbb{P}^{\theta_1}, \mathbb{P}^{\theta_0})/2$, $w > 0$, and $0<\theta_0<\theta<\theta_1<\infty$ such that
$\alpha \Lambda'(\theta)-\Lambda(\alpha)>0$ for $\alpha \in [\theta_0,\theta_1]$,
$\max\{\alpha \Lambda'(\theta-w)-\Lambda(\alpha),~ \alpha\Lambda'(\theta+w)-\Lambda(\alpha)\}<c $ for $\alpha \in [\theta_0,\theta_1]$,
and $F(\alpha)$ has a derivative $F'(\alpha)$, which is positive and continuous for $\theta_0 \leq \alpha \leq \theta_1$.
Since ${\mathbb P}_{1}^\theta \big\{N \geq  (2 \log B)/ K(\mathbb{P}^{\theta}, \mathbb{P}^{\theta_0}) \big\}$
is arbitrarily small when $B$ is large enough, and since for all $C>0,~ {\mathbb E}_1^{\theta}(N|N>C) \leq C+ (2 \log B)/
K(\mathbb{P}^{\theta}, \mathbb{P}^{\theta_0})$ for large enough $B$, it suffices to show that
\begin{eqnarray}\label{39}
&~& (\log B){\mathbb P}_1^{\theta} \bigg\{\displaystyle \max_{n=1,\ldots,(2 \log B)/
K(\mathbb{P}^{\theta_1},\mathbb{P}^{\theta_0})} \int_{\Theta \; \setminus \; [\theta_0, \theta_1]} \\
&~& \times  \displaystyle \sum_{k=1}^{n} \exp \left(\alpha \sum_{i=k}^{n} g(Y_i,Y_{i+1}) -
(n-k+1)\Lambda(\alpha) \right)dF(\alpha) \nonumber \\
&~&~~~~ \geq \displaystyle \frac{4B}{\log B} \bigg\} 
 \longrightarrow_{B \rightarrow \infty }0.\nonumber
\end{eqnarray}

Since the proof of (\ref{39}) follows directly as that in (39) of Pollak \cite{Pollak1987}, we will not repeat it here.

Thus, by Lemmas \ref{61}-\ref{6l}, all conditions of Theorems 3 in Fuh \cite{Fuh2004b} are satisfied,
and so the proof of Theorem 1 is complete.\\ % ~~~~$\Box$~ \\

\def\theequation{7.\arabic{equation}}
\setcounter{equation}{0}
\section{Proof of Theorem \ref{thm1}}\label{sec7}

To prove Theorem \ref{thm1}, we need the following lemmas first. Note that the probability and expected value
are taken under ${\mathbb P}_\omega$ and
${\mathbb E}_\omega$ for $1 \leq \omega < \infty$, we delete $\omega$ for simplicity.
\begin{lemma}\label{L91} Under assumptions of {\rm Theorem \ref{thm1}}.
Let $0<a \leq b < \infty$ satisfy $\Lambda'(a)>\Lambda(b)/b,\\~[a,b]\subset J.$
For any $c > 1$ and probability measure $G$ on $[a,b]$, define
\begin{eqnarray}\label{nbg}
&~& N(c;a,b,G) \\
&=& \inf \bigg\{n | \int_a^b \frac{r(Y_n,\alpha)}{r(Y_0,\alpha)} \exp \{ \alpha {\mathbb S}_n -
n \Lambda (\alpha) \}d G(\alpha) \geq c \bigg\}. \nonumber
\end{eqnarray}
Then there exist constants $0<A,B<\infty$ independent of $c,~G$ such that
\[ {\mathbb E}^{\theta} N(c;a,b,G) \leq  A\log c +B\]
for all $\theta \in [a,b]$ and $c >1$.
\end{lemma}

PROOF. Define $M(\gamma)=\inf \{n| \frac{r(Y_n,\gamma)}{r(Y_0,\gamma)}
\exp\{\gamma {\mathbb S}_n -n\Lambda(\gamma) \} \geq c \}.$
It follows from a simple modification of Lemma 2 of Fuh \cite{Fuh2004a} that
there exists $0<D<\infty$ such that ${\mathbb E}^{\theta} \{{\mathbb S}_{M(\gamma)}-[M(\gamma)\Lambda(\gamma)+ \log c
- \log \frac{r(Y_n,\alpha)}{r(Y_0,\alpha)}]/ \gamma \} \leq D$
uniformly in $\theta \in [a,b],~ \gamma \in [a,b], c>1$. Therefore, by Wald's identity for Markov random walks (cf.
Fuh and Zhang \cite{Fuh2000}) that for all $\theta,~ \gamma \in [a,b]$, there exists a constant $C$
\begin{eqnarray}\label{wald}
&~& {\mathbb E}^{\theta} M(\gamma) \\
&\leq & [(\log c - \log \frac{r(Y_n,\alpha)}{r(Y_0,\alpha)})/\gamma +D]/[\Lambda'(\theta)-\Lambda(\gamma)/\gamma+C] \nonumber \\
 &\leq& [(\log c - \log \frac{r(Y_n,\alpha)}{r(Y_0,\alpha)})/a +D]/[\Lambda'(a)-\Lambda(b)/b+C]. \nonumber
\end{eqnarray}

From $0< r(y,\alpha) < \infty$ for all $y$ via Proposition 1, and $\int_a^b \frac{r(Y_n,\alpha)}{r(Y_0,\alpha)} \exp\{\alpha {\mathbb S}_n -n\Lambda(\alpha)\}dG(\alpha) \geq \min( \frac{r(Y_n,a)}{r(Y_0,a)}
\exp \{a {\mathbb S}_n  - n\Lambda(a)\},~ \frac{r(Y_n,b)}{r(Y_0,b)}\exp \{b{\mathbb S}_n -n\Lambda(b)\})$ it follows that
$N(c;a,b,G) \leq \\ \max(M(a),M(b)) \leq M(a)+M(b).$ This and $(\ref{wald})$ complete the proof of Lemma \ref{L91}.~~
~~~~~~$\Box$

\begin{lemma}\label{L92} For given $0<a \leq b <\infty ,~ [a,b] \subset J,~ \Lambda'(a)>\Lambda(b)/b$, let $G$ be a
probability on $[a,b]$, and denote $F=\gamma F_0 +(1-\gamma) G,$ where $F_0$ is the probability measure wholly concentrated at $\{0\}$ and $\gamma \in (0,1)$.
Consider the optimal stopping problem defined by a prior distribution $F$ on $\theta$
when $Y_0,Y_1,Y_2, \ldots$ are
a sequence of random variables from a state space model satisfying {\rm C1-C4}.
Assume each observation costs $c>0$ if $\theta \neq \theta_0$, zero if
$\theta=\theta_0$, with loss $=1$ for stopping if $\theta=\theta_0$. Then there exists a constant $0<M<\infty$ independent of $c, F$
such that a Bayes procedure {\rm (}with probability one{\rm )} continues sampling
whenever the posterior risk of stopping is at least $Mc$.
\end{lemma}

PROOF.
By making use of a similar procedure as that in pages 2317-2318 of Fuh \cite{Fuh2004b}, a Bayes rule exists.

Let $\infty >Q>A/e$ where A is defined in Lemma \ref{L91} and define $T_{Qc}$ to be the first time $n \leq \infty $ that the posterior
risk of stopping is at most $Qc$. It is sufficient to prove for some $Q<M<\infty$ that the (integrated) risk of $T_{Qc}$
is less than $\gamma$ if $\gamma \geq Mc$. Since the (integrated) risk of any generalized stopping time $T$
is the expected posterior risk of stopping plus $c(1-\gamma) \int_a^b {\mathbb E}^{\theta}T dG(\theta)$,
it is sufficient to prove for some $0<M<\infty$ that $(1-\gamma) \int_a^b {\mathbb E}^{\theta}T_{Qc}dG(\theta)< \gamma/c-Q$ if $\gamma \geq Mc.$

Choose $M>Q$ such that $(1-A/(Qe))M-(B+A/e)>Q$ where $A, B$ are the constants defined by Lemma \ref{L91}.
It is enough to look at $c$ for which $Qc<1$.  Denote
$p_n(\underbar{Y}_n;\theta_0):= p_n(Y_1,\ldots,Y_n;\theta_0)$. 
Note that
\begin{eqnarray}\label{tqc}
&~& T_{Qc} \\
&=& \inf \bigg\{n|Qc \geq \frac{\gamma p_n(\underbar{Y}_n;\theta_0)}{\gamma p_n(\underbar{Y}_n;\theta_0) +(1-\gamma) \int_a^b  p_n(\underbar{Y}_n;\theta) dG(\theta)} \bigg\} \nonumber \\
 &=& \inf \bigg \{ n| \int_a^b \frac{r(W_0,\theta_0)}{r(W_n,\theta_0)} \frac{r(W_n,\alpha)}{r(W_0,\alpha)} \exp \{(\alpha-\theta_0) {\mathbb S}_n \nonumber \\
&~&~~~~~~ - n (\Lambda(\alpha)- \Lambda(\theta_0)\}dG(\alpha)
\geq \displaystyle \frac{\gamma}{1-\gamma} \displaystyle \frac{1-Qc}{Qc} \bigg \} \nonumber \\
 &\leq& \inf \bigg \{ n| \int_a^b  \frac{r(W_0,\theta_0)}{r(W_n,\theta_0)} \frac{r(W_n,\alpha)}{r(W_0,\alpha)} \exp \{(\alpha-\theta_0) {\mathbb S}_n \nonumber \\
&~&~~~~~~ - n (\Lambda(\alpha)- \Lambda(\theta_0)\}dG(\alpha) 
\geq \frac{\gamma}{(1-\gamma)Qc} \bigg \}. \nonumber
\end{eqnarray}

Note that $\sup_{0<\alpha<1}-\alpha (\log \alpha)=1/e$, applying Lemma \ref{L91} to get that if $1>\gamma\geq Mc$
\begin{eqnarray*}
&~& (1-\gamma)\int_a^b {\mathbb E}^{\theta}T_{Qc}dG(\theta)\\
&\leq& (1-\gamma) \left [ A \left ( \log \displaystyle \frac{\gamma}{Qc} + \log \displaystyle \frac{1}{1-\gamma}\right)+B \right]\\
&\leq& \displaystyle \frac{\gamma}{c} \displaystyle \frac{A}{Q} \frac{Qc}{\gamma}\log \displaystyle \frac{\gamma}{Qc}+B+A(1-\gamma)\log \displaystyle \frac{1}{1-\gamma} \\
&\leq& \displaystyle \frac{\gamma}{c} \displaystyle \frac{A}{Qe}+B+\displaystyle \frac{A}{e}\\
& \leq& \displaystyle \frac{\gamma}{c} - \left( 1- \frac{A}{Qe} \right)M+B+\frac{A}{e} \leq
\displaystyle \frac{\gamma}{c}-Q.
\end{eqnarray*}

This completes the proof of Lemma \ref{L92}.~~~~~~~~~~$\Box$\\

PROOF OF THEOREM \ref{thm1}. Without loss of generality, we assume $\theta_0=0$, $0 < a < b$, $(a,b)=J$ and
$\Lambda'(a)>\Lambda(b)/b$. We first show that the right hand side of $(\ref{ao})$ is a lower bound of the
left hand side of $(\ref{ao})$. Consider the Bayesian problem defined in Lemma \ref{L92}
when $\gamma=\frac{1}{2}$ and $dG(\theta)/d\theta = K({\mathbb P}^{\theta}, {\mathbb P}^{\theta_0}) /\int_a^b
K({\mathbb P}^{\alpha}, {\mathbb P}^{\theta_0})  d\alpha$ on $[a,b]$. Let M be the
constant derived in Lemma \ref{L92} and let
$T_{Mc}$ be $T_{Qc}$ for Q=M where $T_{Qc}$ is defined in $(\ref{tqc})$.
$T_{Mc}$ is a mixture stopping rule defined by $G$
and $B =(1-Mc)/(Mc)$. By virtue of Lemma \ref{L92} there exists a Bayes rule
which continues sampling at least as long as $T_{Mc}$. Hence the
Bayes risk is at least the sampling cost of $T_{Mc}$, whence for any stopping rule T
\[ {\mathbb P}^{\theta_0}(T<\infty) + c\int_a^b {\mathbb E}^{\theta} T dG(\theta)
\geq c\int_a^b {\mathbb E}^{\theta}T_{Mc}dG(\theta).\]
Thus if ${\mathbb P}^{\theta_0}(T<\infty) \leq 1/\varepsilon = Mc/(1-Mc),$ then
\begin{eqnarray}\label{etc}
\int_a^b {\mathbb E}^{\theta} T dG(\theta) \geq \int_a^b {\mathbb E}^{\theta} T_{Mc} dG(\theta)-M/(1-Mc).
\end{eqnarray}
There exist $a_1,b_1$ such that $0<a_1<a<b<b_1<\infty$ and $\Lambda'(a_1)>\Lambda(b_1)/b_1$.
Define $\Lambda=\inf\{n|\int_{a_1}^{b_1} \exp \{\theta {\mathbb S}_n-n \Lambda(\theta)\}K({\mathbb P}^{\theta},
{\mathbb P}^{\theta_0})  d\theta/ \int_a^b K({\mathbb P}^{\theta}, {\mathbb P}^{\theta_0})
d\theta \geq B \}$. By definition,
$T_{Mc} \geq \Lambda$. $\Lambda$ is a mixture stopping rule defined by $dF(\theta)/d\theta
=K({\mathbb P}^{\theta}, {\mathbb P}^{\theta_0})
/\int_{a_1}^{b_1}K({\mathbb P}^{\alpha}, {\mathbb P}^{\theta_0})  d\alpha$ on $[a_1,b_1]$
and $B'= B \int_a^b K({\mathbb P}^{\theta}, {\mathbb P}^{\theta_0})
d\theta/\int_{a_1}^{b_1}K({\mathbb P}^{\theta},{\mathbb P}^{\theta_0})  d\theta$.
Thus by Theorem 1
\begin{eqnarray}\label{etmc}
&~& {\mathbb E}^{\theta} T_{Mc} \geq {\mathbb E}^{\theta} \Lambda \\
&=& \frac{1}{2K({\mathbb P}^{\theta},{\mathbb P}^{\theta_0}) }[2\log B' + \log \log B']  +O_{\theta}(1), \nonumber
\end{eqnarray}
where $ \lim \sup_{\varepsilon \rightarrow \infty} \sup_{a \leq \theta \leq b}|O_{\theta}(1)| \leq \infty$.
Combining $(\ref{etc})$ and $(\ref{etmc})$, and replacing $B'$ by $B$ yields
\begin{eqnarray*}
&~& \int_a^b {\mathbb E}^{\theta} T dG(\theta) \\
&\geq& \int_a^b [2 \log B +\log \log B +O(1)]d\theta/2(\int_a^b K({\mathbb P}^{\theta}, {\mathbb P}^{\theta_0})  d\theta).\end{eqnarray*}
Hence by definition of $G$, we have
\[ \int_a^b[2 K({\mathbb P}^{\theta}, {\mathbb P}^{\theta_0})  {\mathbb E}^{\theta} T -
(2\log B +\log \log B )]d\theta \geq O(1) \]
for all $T$ satisfying ${\mathbb P}^{\theta_0}\{T<\infty\} \leq 1/ B$.

To show that the equality is attained by the weighted SRP rule. By Theorem 1, we need only to show
that for the weighted SRP detection rule {\rm (\ref{srp})} satisfies
${\mathbb P}^{\theta_0} \{N_b < \infty\} \leq 1/c,$
%~~~{\rm and~~~}{\mathbb P}^{\theta} \{N_b < \infty\}=1~~~{\rm for~all~}\theta \in J,
for any $c >1$.

Recall that ${\mathbb P}^{\theta}(y,dz)$ defined in (\ref{et}), and denote
${\mathbb Q}(w,dw'):=\int_{\theta \in J} {\mathbb P}^{\theta}(w,dw') dF(\theta)$.
Then it is easy to see that ${\mathbb Q}(w,\cdot)$ is a transition kernel.
By definition of $N_b$, we have
\begin{eqnarray}
&~& {\mathbb P}^{\theta_0} \{N_b < \infty\} \\
&=& \int_{\{N_b < \infty\}} \frac{1}{LR_n(F)} d{\mathbb Q}
\leq \frac{1}{c} {\mathbb Q}\{N_b < \infty\}  \leq \frac{1}{c}. \nonumber
\end{eqnarray}
This establishes the desired property, and thus completing the proof of $(\ref{ao})$.~~~~~~~~~~~~$\Box$\\

\def\theequation{A.\arabic{equation}}
\setcounter{equation}{0}
\centerline{APPENDIX}
%{\bf Proofs of Propositions \ref{A1}.}

We give a proof of Proposition \ref{A1} which also corrects notations error in Lemma 8 of Fuh \cite{Fuh2004b},
in the setting of hidden Markov models.

PROOF OF PROPOSITION \ref{A1}. Note that
\begin{eqnarray*}
&~& H_{n+1}(y,w) \\
&=& \mathbb{P}_{\infty} \{ R_{n+1,p} \leq y|  N_{q,b} > n+1, W_{n+1} =w\}  \\
&=& \int_{w'\in {\cal W}} \int_0^{B(w')} \mathbb{P}_{\infty} \{ R_{n+1,p} \leq y, W_{n} \in dw', \\
&~&~~~~~~~~~~~R_{n,p}\in dt|  N_{q,b} > n +1, W_{n+1}=w\} \\
&=& \int_{w'\in {\cal W}} \int_0^{B(w')} \mathbb{P}_{\infty} \{ R_{n+1,p} \leq y| W_{n} = w', R_{n,p}=t, \\
&~&~~~~~~~~~~~    N_{q,b} > n +1, W_{n+1}=w\} \\
&~&~\times \mathbb{P}_{\infty} \{ R_{n,p} \in dt, W_{n} \in dw'|
    N_{q,b} > n +1, W_{n+1}=w\} \\
&=& \int_{w'\in {\cal W}} \int_0^{B(w')} \rho (t,y,w) \mathbb{P}_{\infty} \{ R_{n,p} \in dt, W_{n} \in dw'| \\
&~&~~~~~~~~~~~ N_{q,b} > n +1, W_{n+1}=w\}.
\end{eqnarray*}
Since
\begin{eqnarray*}
&~& \mathbb{P}_{\infty} \{ R_{n,p} \in dt, W_{n} \in dw' |  N_{q,b} > n+1, W_{n+1} =w\}  \\
&=& \mathbb{P}_{\infty} \{ R_{n,p} \in dt,W_{n} \in dw' | N_{q,b} > n, \\
&~&~~~~~~ N_{q,b} > n+1, W_{n+1} =w\} \\
&=& \bigg(\mathbb{P}_{\infty} \{ N_{q,b} > n+1,W_{n+1} \in dy| R_{n,p}=t, N_{q,b} > n, \\
&~& W_n=w'\} \bigg) \bigg/ \bigg( \int_{w,w'\in {\cal W}} \int_0^{B(w')}  \mathbb{P}_{\infty} \{ N_{q,b} > n+1,\\
&~&      W_{n+1} \in dw| R_{n,p}=t,   N_{q,b} > n, W_n=w'\} \\
&~&~~~ \times \frac{\mathbb{P}_{\infty} \{ R_{n,p} \in dt|  N_{q,b} > n, W_n =w'\}}
      { \mathbb{P}_{\infty} \{R_{n,p} \in dt |  N_{q,b} > n, W_n =w'\}} \\
&~&~~~ \times \frac{ \mathbb{P}_{\infty} \{ N_{q,b} > n| W_n =w'\} \mathbb{P}(w',dw)}{ \mathbb{P}_{\infty} \{ N_{q,b} > n| W_n =w'\} \mathbb{P}(w',dw)} \bigg)\\
&=& \frac{\zeta(t,w',w)dH_n(t,w')\mathbb{P}(w',dw)}{\int_{w,w'\in {\cal W}} \int_0^{B(w')}
    \zeta(t,w',w)dH_n(t,w')\mathbb{P}(w',dw)}.
\end{eqnarray*}
It follows that
\begin{eqnarray*}
&~&H_{n+1}(y,w) \\
&=& \frac{ \int_{w' \in {\cal W}} \int_0^{B(w')} \rho(t,y,w) \zeta(t,w',w) dH_n(t,w')\mathbb{P}(w',dw) } {Q(H_n)} \\
&=& T_BH_n(y,w).
\end{eqnarray*}

The existence of the fixed point follows the same argument as that of Lemma 11 in Pollak \cite{Pollak1985}. $\hfill \Box$~\\

{\bf Acknowledgement}: 
The author is grateful for Professor Alex Tatakovsky's valuable comments.

\end{document}